\documentclass[11pt,leqno]{amsart}

\usepackage[colorlinks=false,hidelinks]{hyperref} 
\usepackage[T1,T5]{fontenc} 
\usepackage[vietnamese,english]{babel}
\usepackage[colorlinks=false,hidelinks]{hyperref} 
\usepackage[T1,T5]{fontenc} 
\usepackage[vietnamese,english]{babel}
\usepackage{mathtools}
\usepackage{amssymb, amsmath, amsthm, amsfonts}
\usepackage[shortlabels]{enumitem}
\usepackage{graphicx}
\usepackage{orcidlink}
\usepackage{quiver}
\usetikzlibrary{babel}
\usepackage[initials]{amsrefs}
\usepackage{quiver}
\usepackage{lmodern}
\usepackage[all,cmtip]{xy}
\usepackage{longtable}
\usepackage[english,ruled,vlined, algosection]{algorithm2e}
\usepackage{caption}

\hypersetup{
	pdflang={en-US},
	pdftitle={Classification and structure of generalized Legendrian racks},
	pdfauthor={Lực Ta},
	pdfsubject={Mathematics},
	pdfkeywords={Center, classification, generalized Legendrian rack, quandle, rack, tensor product	}
}

\usetikzlibrary{automata,positioning}
\DeclareMathOperator{\id}{id}
\newcommand{\R}{\mathbb{R}}
\newcommand{\Z}{\mathbb{Z}}

\newtheorem{thm}{Theorem}[section]
\newtheorem{lemma}[thm]{Lemma}
\newtheorem{cor}[thm]{Corollary}
\newtheorem{prop}[thm]{Proposition}

\theoremstyle{definition}
\newtheorem{definition}[thm]{Definition}
\newtheorem{example}[thm]{Example}
\newtheorem*{competing}{Competing Interests}
\newtheorem*{ack}{Acknowledgments}
\theoremstyle{remark}
\newtheorem{remark}[thm]{Remark}

\renewcommand{\phi}{\varphi}

\newcommand{\Set}{\mathsf{Set}}
\newcommand{\Grp}{\mathsf{Grp}}
\newcommand{\Rc}{\mathsf{Rack}}
\newcommand{\Qn}{\mathsf{Qnd}}
\newcommand{\glr}{\mathsf{GLR}}
\renewcommand{\glq}{\mathsf{GLQ}}
\newcommand{\gla}{{\mathsf{GLR}_\mathrm{med}}}

\newcommand{\F}{F}
\newcommand{\G}{G}
\newcommand{\Inn}{\operatorname{Inn}}
\newcommand{\Conj}{\operatorname{Conj}}

\newcommand{\inv}{^{-1}}
\newcommand{\tr}{\triangleright}
\renewcommand{\u}{\mathtt{u}}
\newcommand{\downcusp}{\mathtt{d}}
\newcommand{\ud}{\mathtt{ud}}
\newcommand{\du}{\mathtt{du}}

\newcommand{\GLR}{\mathcal{G}}
\DeclareMathOperator{\Hom}{Hom}
\DeclareMathOperator{\Aut}{Aut}
\DeclareMathOperator{\Free}{Free}

\newcommand{\med}{_{\mathrm{med}}}
\newcommand{\perm}{_{\mathrm{perm}}}
\newcommand{\qnd}{_{\mathrm{qnd}}}
\newcommand{\GG}{\mathcal{G}}

\newcommand{\Fs}{{\widetilde{s}}}
\newcommand{\Ft}{{\widetilde{t}}}
\newcommand{\Gs}{{\hat{s}}}
\newcommand{\Gt}{{\hat{t}}}

\DeclareMathOperator{\Forg}{For}
\newcommand{\blr}{\mathsf{BLR}}

\begin{document}
	
	\captionsetup[longtable]{labelfont=sc}
	
	\title[Classification and structure of GL-racks]{Classification and structure\\of generalized Legendrian racks}
	\author{L\d\uhorn c Ta}	
	
	\address{Department of Mathematics, University of Pittsburgh, Pittsburgh, Pennsylvania 15260}
	\email{ldt37@pitt.edu}
	
	\subjclass[2020]{Primary 20N02; Secondary 08A35, 18B99, 20B25, 57K12}
	
	\keywords{Center, classification, generalized Legendrian rack, quandle, rack, tensor product}
	
	\begin{abstract}
	This paper develops the structure theory of generalized Legendrian racks, which are algebraic structures used to study Legendrian knots. First, we give a group-theoretic characterization of GL-racks. As applications, we classify GL-racks up to order 8 up to isomorphism via an exhaustive search algorithm, and we classify certain infinite families of GL-racks. Then we compute the center of the category of GL-racks and construct an equivalence of categories between racks and GL-quandles. Finally, we study tensor products of racks and GL-racks coming from universal algebra.
	\end{abstract}
	\maketitle
	
	\section {Introduction}
\emph{Generalized Legendrian racks}, also called \emph{GL-racks} or \emph{bi-Legendrian racks}, are algebraic structures used to distinguish Legendrian links in contact three-space. GL-racks can be traced back to algebraic structures called \emph{kei}, which Takasaki \cite{takasaki} introduced in 1942 to study symmetric spaces; \emph{quandles}, which Joyce \cite{joyce} and Matveev \cite{matveev} independently introduced in 1982 to study links in $\R^3$ and $S^3$; and \emph{racks}, which Fenn and Rourke \cite{fenn} introduced in 1992 to study framed links in $3$-manifolds.
Kei, quandles, and racks have enjoyed significant study as link invariants in geometric topology and in their own rights in quantum algebra and group theory.

More recently, various authors have introduced variants of racks suitable for studying Legendrian links. In 2017, Kulkarni and Prathamesh \cite{original} introduced rack invariants of Legendrian knots. In 2021, Ceniceros, Elhamdadi, and Nelson \cite{ceniceros} refined these invariants by introducing \emph{Legendrian racks}. 
In 2023, Karmakar, Saraf, and Singh \cite{karmakar} and Kimura \cite{bi} independently strengthened these constructions by introducing GL-racks, which are racks equipped with additional structure. In 2026, Cheng and He \citelist{\cite{cheng}\cite{cheng2}} studied the efficacy of GL-racks as Legendrian knot invariants, and Karmakar, Saraf, and Singh \cite{karmakar2} studied Beck modules over GL-racks. The author and Wood \cite{wood} studied \emph{4-Legendrian racks}, of which GL-racks form a special class. The author \cite{ta} also found relationships between involutory GL-racks and similar algebraic structures used to study surface-knots.

The purpose of this paper is to develop the structure theory of GL-racks in hopes of strengthening their applications in Legendrian knot theory. In particular, we study group-theoretic, categorical, and universal-algebraic aspects of GL-racks, and we classify GL-racks up to order $8$.

\subsection{Organization of the paper}
	In Section \ref{sec:racks}, we define racks and quandles, consider several examples, and discuss a canonical rack automorphism $\pi$ that plays a fundamental role in the theory. 

In Section \ref{sec:gl-racks}, we give a simplified definition of GL-racks, show its equivalence to the definition in the literature, and discuss examples and universal-algebraic aspects of GL-racks. In particular, we show how to construct free GL-racks on arbitrary generating sets.

In Section \ref{subsec:classn}, we study group-theoretic aspects of GL-racks. We determine all GL-structures on a given rack, answering a question from a previous version of \cite{karmakar}. As an application, we classify GL-racks up to order $8$ up to isomorphism via a computer algorithm; we provide the program and data in a GitHub repository \cite{my-code}. Then we classify all GL-structures on permutation racks, certain conjugation quandles, and certain Takasaki kei, including all dihedral quandles. We also compute automorphism groups of dihedral GL-quandles. 

	In Section \ref{sec:cat-equiv}, we study categorical aspects of GL-racks. First, we compute the centers of the categories of GL-racks, GL-quandles, Legendrian racks, and Legendrian quandles. Then we construct an equivalence of categories between the categories of racks and GL-quandles. This equivalence restricts to an equivalence between the full subcategories of medial racks and medial GL-quandles. Our construction corresponds to an isomorphism of algebraic theories.

	In Section \ref{sec:tensors}, we study tensor products of racks and GL-racks coming from universal algebra. We show that the categories of racks and GL-racks have tensor units. This is unusual for noncommutative algebraic theories; for example, the tensor product in the category of groups does not have a tensor unit. The induced symmetric monoidal structure on medial racks is closed, and similarly for medial GL-racks.

In Section \ref{sec:questions}, we propose questions for further research.

\begin{ack}
		An early version of this paper was written in partial fulfillment of the senior requirements for the bachelor's degree in mathematics at Yale University. I thank Sam Raskin, my bachelor's thesis advisor, for mentorship and support throughout the research and writing processes. 
		I thank the organizers of the UnKnot V conference in 2024, where I met Jose Ceniceros and Peyton Wood, whom I thank for introducing me to various aspects of the theory. I thank Noam Scully for helpful discussions. I thank Zhiyi He for identifying errors in an earlier draft; these errors were entirely my own. I thank Qiaochu Yuan for a terminological correction and an anonymous Stack Exchange user for referencing me to \cite{alg-theories}. Finally, I thank the anonymous referees for helpful suggestions, particularly for suggesting Theorem \ref{thm:free-arbitrary}.
\end{ack}
	
	\section{Racks}\label{sec:racks}
	\subsection{Racks}
	
	In this subsection, we define racks and discuss several major examples.
	
	\subsubsection{Notation}
	Given a set $X$, we denote the permutation group of $X$ by $S_X$. 
	In the case that $X=\{1,2,\dots,n\}$, we denote the symmetric group on $n$ letters by $S_n$.
	We also denote the composition of functions $\phi:X\to Y$ and $\psi:Y\to Z$ by $\psi\phi$.
	
	While racks and quandles are often defined as sets $X$ with a right-distributive nonassociative binary operation $\tr:X\times X\to X$ satisfying certain axioms, they may also be characterized in terms of permutations $s_x\in S_X$ assigned to each element $x\in X$; cf.\ \citelist{\cite{survey}*{Def.\ 2.1}}. 
	One may translate between the two conventions via the formula \[x\tr y=s_y(x).\] In this paper, we adopt the convention using permutations due to its convenience for abstract proofs and exhaustive search algorithms.
	
	\subsubsection{Racks and quandles}
	Having established our notation, we define racks and quandles. 
	Although we provide all relevant definitions and preliminaries, we also refer the reader to \citelist{\cite{intro}\cite{quandlebook}} for accessible introductions to quandles, \citelist{\cite{quandlebook}\cite{book}} for references on racks as they concern low-dimensional topology, and \cite{survey} for a survey of modern algebraic literature on racks.
	
	\begin{definition}
		Let $X$ be a set, let $s:X\to S_X$ be a map, and write $s_x:=s(x)$ for all elements $x\in X$. We call the pair $(X,s)$ a \emph{rack} if \[s_xs_y=s_{s_x(y)}s_x\] for all $x,y\in X$, in which case we call $s$ a \emph{rack structure} on $X$. 
		If in addition $s_x(x)=x$ for all $x\in X$, then we say that $(X,s)$ is a \emph{quandle}.
		We also say that $|X|$ is the \emph{order} of $(X,s)$.
		Finally, if $Y\subseteq X$ and $s_y^{\pm 1}(z)\in Y$ for all $y,z\in Y$, then we say that $(Y,s|_Y)$ is a \emph{subrack} of $(X,s)$.
	\end{definition}
	
	\begin{definition}\phantom{-}
		\begin{itemize}
			\item Given two racks $R:=(X,s)$ and $(Y,t)$, we say that a map $\phi:X\to Y$ is a \emph{rack homomorphism} if \[\phi s_x=t_{\phi(x)}\phi\] for all $x\in X$. Let $\Rc$ denote the category of racks and rack homomorphisms, and let $\Qn$ be the full subcategory of $\Rc$ whose objects are quandles.
			\item A \emph{rack isomorphism} is a bijective rack homomorphism. The \emph{automorphism group} of $R$ is denoted by $\Aut R$.
			\item The \emph{inner automorphism group} of $R$ is the subgroup $\operatorname{Inn}R:=\langle s_x\mid x\in X\rangle$ of $\Aut R$.
			\item Equip $X\times X$ with the product rack structure. If the natural map $X\to X\times X$ is a rack homomorphism, then we say that $R$ is \emph{medial}. Let $\Rc\med$ be the full subcategory of $\Rc$ whose objects are medial.
		\end{itemize}
\end{definition}
	
	\begin{example}[\cite{quandlebook}*{Ex.\ 99}]
		Let $X$ be a set, and fix $\sigma\in S_X$. Define $s:X\rightarrow S_X$ by $x\mapsto \sigma$, so that $s_x(y)=\sigma(y)$ for all $x,y\in X$. Then $(X,\sigma)\perm:=(X,s)$ is a medial rack called a \emph{permutation rack} or \emph{constant action rack}. (Our notation $(X,\sigma)\perm$ is nonstandard, unlike the notation in the next three examples.) Note that $(X,\sigma)\perm$ is a quandle if and only if $\sigma=\id_X$, in which case we call $(X,\id_X)\perm$ a \emph{trivial quandle}.
	\end{example}
	
	\begin{example}[\cite{quandlebook}*{Ex.\ 54}]\label{ex:takasaki-kei}
		Let $A$ be an abelian additive group. Define $s:A\to S_A$ by $b\mapsto s_b$ with $s_b(a):=2b-a$ for all elements $a,b\in A$. Then $T(A):=(A,s)$ is a medial quandle called a \emph{Takasaki kei}. Takasaki kei are the earliest examples of racks in the literature; Takasaki \cite{takasaki} introduced them in 1943 to study symmetric spaces.
		
		In particular, if $A=\Z/n\Z$, then we call $R_n:=T(\Z/n\Z)$ the \emph{dihedral quandle} of order $n$.
	\end{example}
	
	\begin{example}[\cite{book}*{Ex.\ 2.13}]\label{ex:conj}
	Let $X$ be a union of conjugacy classes in a group $G$, and 
	define $c^G:X\rightarrow S_X$ by sending $x$ to the conjugation map \(c^G_x:=[y\mapsto xyx\inv]\).
		Then $\Conj X:=(X,c^G)$ is a quandle called a \emph{conjugation quandle} or \emph{conjugacy quandle}.
\end{example}
	
	\begin{definition}[\cite{dual}]
		Given a rack $R=(X,s)$, define $s':X\to S_X$ by $x\mapsto s_x\inv$. Then $R^{\operatorname{op}}:=(X,s')$ is a rack called the \emph{dual rack} of $R$.
	\end{definition}
	
	The following is straightforward.
	
	\begin{prop}\label{prop:aut-inn-dual}
		For all racks $R=(X,s)$, we have $\Aut R=\Aut R^{\operatorname{op}}$ and $\Inn R=\Inn R^{\operatorname{op}}$. Moreover, $R$ is medial if and only if $R^{\operatorname{op}}$ is medial.
	\end{prop}
	
	\subsection{The canonical automorphism of a rack}\label{sec:canonical-aut}
	In this subsection, we discuss a collection $\Pi$ of canonical rack automorphisms $\pi_R\in\Aut R$ that generates the center of the category $\Rc$.

\subsubsection{Definition of $\pi_R$}
Given a rack $R=(X,s)$, define $\pi_R:X\to X$ by ${x\mapsto s_x(x)}$; see \cite{center}*{Prop.\ 2.5}. This function is sometimes called the \emph{kink map} because of its relationship to kinks in framed knot diagrams; see \cite{quandlebook}*{p.\ 149}.

Evidently, $R$ is a quandle if and only if $\pi_R=\id_X$. Thus, we can loosely think of $\pi_R$ as measuring the failure of $R$ to be a quandle. Whenever there is no ambiguity, we will suppress the subscript and only write $\pi:=\pi_R$.

\subsubsection{Categorical centers}
Recall that the \emph{center} of a category $\mathcal{C}$ is the commutative monoid $Z(\mathcal{C})$ of natural endomorphisms of the identity functor $\mathbf{1}_\mathcal{C}$. 
Concretely, $\eta\in Z(\mathcal{C})$ if and only if, for all objects $R,S$ and morphisms $\phi:R\to S$ in $\mathcal{C}$, the component $\eta_R$ is an endomorphism of $R$, and $\eta_S\phi=\phi\eta_R$. 

For example, if $A\text{-}\mathsf{mod}$ denotes the category of modules over a ring $A$, then the categorical center $Z(A\text{-}\mathsf{mod})$ is isomorphic to the ring-theoretic center $Z(A)$ of $A$.

\subsubsection{Properties of $\pi$}
Let $\Pi$ denote the collection of canonical automorphisms $\pi_R$ for all racks $R$. 
Szymik \cite{center} proved the claims in the following proposition using the binary operation $\tr$. We offer a new proof of claim \ref{iii} and rewrite Szymik's proofs of the other claims in terms of permutations.

\begin{prop}\label{lem:theta-inv}
	For all racks $R=(X,s)$ and all integers $k\in\Z$, we have the following:
	\begin{enumerate}[({A}1)]
		\item\label{a2} $\pi:X\to X$ is a bijection with inverse $\pi\inv$ defined by $x\mapsto s_x\inv(x)$.
		\item\label{a1} $\pi^{\pm 1}_S\phi=\phi\pi^{\pm 1}_R$ for all racks $S=(Y,t)$ and rack homomorphisms $\phi\in\Hom_\Rc(R,S)$.
		\item\label{iii} $s_{\pi^k(x)}=s_x$ for all $x\in X$.
		\item\label{a3} $\Pi^k\in Z(\Rc)$.
	\end{enumerate}
\end{prop}

\begin{proof}
	Define $\pi\inv$ as in claim \ref{a2}, and fix $x\in X$. We deduce from Proposition \ref{prop:aut-inn-dual} that
	\[
	\pi\inv\pi(x)=\pi\inv s_x(x)=s\inv_{s_x(x)}s_x(x)=s_xs\inv_x(x)=x,
	\]
	as desired. Dually, $\pi\pi\inv(x)=x$, which proves claim \ref{a2}. 
	To prove claim \ref{a1}, observe that
	\[
	\pi_S\phi(x)=t_{\phi(x)}\phi(x)=\phi s_x(x)=\phi\pi_R(x).
	\]
	Since $\pi_S$ and $\pi_R$ are bijections, we obtain claim \ref{a1}.
	
	To prove claim \ref{iii}, note that the rack structure $s:X\to S_X$ is a rack homomorphism from $R$ to $\Conj S_X$. By claim \ref{a1} and the fact that $\Conj S_X$ is a quandle,
	\begin{equation*}
		s_{\pi^k_{R} (x)}=s\pi^k_R(x)=\pi^k_{\Conj S_X}s(x)=s(x)=s_x
	\end{equation*}
	for all integers $k\in\Z$, which proves claim \ref{iii}. 
	Now, claims \ref{a1} and \ref{iii} yield
	\[
	s_{\pi(x)}\pi=\pi s_{\pi(x)}=\pi s_x,
	\]
	so $\pi$ is a rack endomorphism. Combined with claims \ref{a2} and \ref{a1}, this proves claim \ref{a3}.
\end{proof}

Hence, $\Pi$ generates a cyclic subgroup of $Z(\Rc)$. In fact, this subgroup \emph{is} $Z(\Rc)$; see \cite{center}*{Thm.~5.4}.

The following identities will be useful later on.

\begin{prop}\label{prop:rack-center}
	Let $R=(X,s)$ be a rack. Then the following hold for all $k\in\Z$ and $x\in X$:
	\begin{enumerate}[({B}1)]
		\item\label{i} $\pi^k s_x=s_x\pi^k$.
		\item\label{ii} $\pi^k(x)=s^k_x(x)$.
	\end{enumerate}
\end{prop}

\begin{proof}
	Claim \ref{i} is immediate from the inclusions $s_x\in\Aut R$ and $\Pi^{\pm 1}\in Z(\Rc)$.
	
	To prove claim \ref{ii}, we induct on $k$. The base case $k=0$ is trivial. For $k>0$, we have
	\[
	\pi^k(x)=\pi^{k-1}\pi(x)=\pi^{k-1} s_x(x)=s_x\pi^{k-1}(x)=s_xs^{k-1}_x(x)=s^k_x(x)
	\]
	as desired; the second equality uses the definition of $\pi$, while the third equality follows from claim \ref{i}. Thanks to part \ref{a2} of Proposition \ref{lem:theta-inv}, a similar argument proves claim \ref{ii} for $k<0$.
\end{proof}
	
	\subsection{Racks as an algebraic theory}\label{sec:rack-thy}
	
	Recall that racks form an algebraic theory with two binary operations $s_-(-)$ and $s\inv_-(-)$. In this subsection, we discuss universal-algebraic free objects in $\Rc$ and rack structures on certain hom-sets of racks. While we provide the relevant results here, we also direct the reader to \cite{alg-theories} for a reference on the underlying universal algebra.
	
	\subsubsection{Examples of quotient racks} Since $\Rc$ is the category of models of the algebraic theory of racks in $\Set$, \cite{alg-theories}*{Thm.\ 3.4.5} states that $\Rc$ is complete and cocomplete. Thus, we can consider quotients of racks by congruences; see \cite{alg-theories}*{Lem.\ 3.5.1} and cf.\ \cite{dual}. This allows us to state the following definitions, which are actually left adjoints of the inclusion functors from $\Qn$ and $\Rc\med$ into $\Rc$.
	
	\begin{definition}[\cite{fenn}*{Sec.\ 2}, \cite{joyce}*{Sec.\ 10}]
		Given a rack $R=(X,s)$, the \emph{associated quandle} of $R$ is the quandle $R\qnd:=R/{\sim}$, where $\sim$ is the smallest congruence on $X$ such that $s_x(x)=x$ for all $x\in X$. Similarly, the \emph{medialization} or \emph{abelianization} of $R$ is the medial rack $R\med:=R/{\sim}$, where $\sim$ is the smallest congruence on $X$ such that 
		\begin{equation}\label{eq:medial-2}
			{s_{s_x(z)}s_y(a)=s_{s_x(y)}s_z(a)}
		\end{equation}
		 for all $x,y,z,a\in X$. 
	\end{definition}
	
	\subsubsection{Free racks}
	In analogy to free groups or free modules, the \emph{free rack} $\Free_\Rc (X)$ on a set $X$ exists and is uniquely (up to isomorphism) characterized by a universal property; see, for example, \cite{alg-theories}*{Cor.\ 3.7.8}. Namely, there exists a set map $i:X\to\Free_\Rc(X)$ such that for all racks $R=(Y,t)$ and all set maps $\phi_0:X\to Y$, there exists a unique rack homomorphism ${\phi\in\Hom_\Rc(\Free_\Rc(X),R)}$ such that $\phi i=\phi_0$. That is, the following diagram commutes:
	\[\begin{tikzcd}
		X & {\operatorname{Free}_{\Rc}(X)} \\
		& Y
		\arrow["i", from=1-1, to=1-2]
		\arrow["{\phi_0}"', from=1-1, to=2-2]
		\arrow["\phi", dashed, from=1-2, to=2-2]
	\end{tikzcd}\]
	For an explicit set-theoretic construction of free racks, see \cite{free-z}*{Prop.\ 1.3}. We can also define the \emph{free quandle} on $X$ to be $(\Free_\Rc(X))\qnd$; cf.\ \cite{fenn}*{Sec.\ 7.1}.
	
	\begin{example}[\cite{free-z}*{Ex.\ 1.6}]\label{ex:free-one}
		The free rack $F$ on one generator is canonically isomorphic to the permutation rack $(\Z,\sigma)\perm$, where $\sigma(k)=k+1$ for all $k\in\Z$. 
		Under this identification, $\pi_F=\sigma$. 
		By contrast, the free quandle on one generator is the trivial quandle with one element.
		
		For fun, we encourage the reader to prove that $F\cong (\Z,\sigma)\perm$. To do this, appeal to the universal property of $\Free_\Rc(\{x\})$ with $i(x):=0$ and $\phi(k):=\pi_R^k\phi_0(x)$.
	\end{example}
	
	\subsubsection{Commutative algebraic theories} 
	Medial racks are the largest commutative subtheory of racks. Recall that an algebraic theory $\mathcal{T}$ is called \emph{commutative} if, for all $\mathcal{T}$-models $A$, every $n$-ary operation $\alpha:\mathcal{T}^n\to \mathcal{T}$ defines a homomorphism $A^n\to A$; see \cite{alg-theories}*{Thm.\ 3.10.3}. 

Let $\mathcal{T}_{\operatorname{ab}}$ be the largest commutative subtheory of an algebraic theory $\mathcal{T}$. For example, if $\mathcal{T}$ is the algebraic theory of groups, then $\mathcal{T}_{\operatorname{ab}}$ is the algebraic theory of abelian groups. 
If $\mathcal{C}$ is the category of $\mathcal{T}$-models in $\Set$, let $\mathcal{C}_{\operatorname{ab}}$ be the full subcategory of $\mathcal{C}$ whose objects are $\mathcal{T}_{\operatorname{ab}}$-models. The following strengthens part of \cite{alg-theories}*{Thm.\ 3.10.3}; we could not find a reference for this precise statement.
	
	\begin{prop}\label{prop:hom-alg-thy}
		For all $\mathcal{T}$-models $X$ in $\mathcal{C}$ and $\mathcal{T}_{\operatorname{ab}}$-models $Y$ in $\mathcal{C}_{\operatorname{ab}}$, the set $H:=\Hom_\mathcal{C}(X,Y)$ has a canonical $\mathcal{T}_{\operatorname{ab}}$-model structure defined by \[\alpha(f_1,\dots,f_n)(x):=\alpha(f_1(x),\dots,f_n(x))\] for all $n$-ary operations $\alpha:\mathcal{T}^n\to\mathcal{T}$, $\mathcal{T}$-model homomorphisms $f_1,\dots,f_n\in H$, and elements $x\in X$.
	\end{prop}
	
	\begin{proof}
		We refer the reader to the proof of the statement in \cite{alg-theories}*{Thm.\ 3.10.3} that (1) implies (3). This proof uses the commutativity of $n$-ary operations in $Y$ but not in $X$, so the implication still holds with the weakened assumption that $X$ is not necessarily a $\mathcal{T}_{\operatorname{ab}}$-model.
	\end{proof}
	
	\subsubsection{Mediality and Hom racks}
	If we take $\mathcal{T}$ to be the algebraic theory of racks (resp.\ quandles), then $\mathcal{T}_{\operatorname{ab}}$ is the algebraic theory of medial racks (resp.\ medial quandles). That is, if $\mathcal{C}=\Rc$, then $\mathcal{C}_{\operatorname{ab}}=\Rc\med$; this follows directly from Proposition \ref{prop:aut-inn-dual}. 
	
	Hence, Proposition \ref{prop:hom-alg-thy} recovers the following result of Grøsfjeld \cite{rack-roll}*{Prop.\ 3.3} in 2021, which in turn generalized a result of Crans and Nelson \cite{Hom}*{Thm.\ 3} in 2014.
	
	\begin{cor}\label{cor:hom-racks}
		Let $R=(X,s)$ be a rack, and let $M=(Y,t)$ be a medial rack. Then the hom-set $H:=\Hom_\Rc(R,M)$ has a canonical medial rack structure $\Tilde{t}:H\to S_H$ defined by \[\Tilde{t}_g(f)(x):=t_{g(x)}f(x)\] for all $f,g\in H$ and $x\in X$. Moreover, if $R$ is a quandle or $M$ is a quandle, then $H$ is also a quandle.
	\end{cor}
	
	\begin{proof}
		The only statement that does not directly follow from Proposition \ref{prop:hom-alg-thy} is the final sentence, which is straightforwardly verified.
	\end{proof}
	
	\section{Generalized Legendrian racks}\label{sec:gl-racks}
	\subsection{GL-racks}
	In this subsection, we define \emph{generalized Legendrian racks} (also called \emph{GL-racks}), which Karmakar, Saraf, and Singh \cite{karmakar} and Kimura \cite{bi} introduced independently in 2023.
	Once again, we express our definition in terms of permutations.
	
	\begin{definition}\label{def:gl-rack}
		Given a rack $R=(X,s)$, a \emph{generalized Legendrian structure} or \emph{GL-structure} on $R$ is a rack automorphism $\u\in\Aut R$ such that $\u s_x=s_x\u$ for all $x\in X$.
		We call the pair $(R,\u)$ a \emph{generalized Legendrian rack} or \emph{GL-rack}. If in addition $R$ is a quandle or a medial rack, then we also call $(R,\u)$ a \emph{GL-quandle} or a \emph{medial GL-rack}, respectively.
	\end{definition}
	
	\begin{definition}\label{def:gl-hom}
		A \emph{GL-rack homomorphism} between two GL-racks $(R_1,\u_1)$ and $(R_2,\u_2)$ is a rack homomorphism $\phi\in\Hom_\Rc(R_1,R_2)$ satisfying $\phi\u_1=\u_2\phi$. We denote the category of GL-racks and their homomorphisms by $\glr$. Finally, let $\glq$, $\glr\med$, and $\glq\med$ be the full subcategories of $\glr$ whose objects are GL-quandles, medial GL-racks, and medial GL-quandles, respectively.
	\end{definition}
	
	\begin{remark}\label{rmk:virtual}
		\emph{Virtual racks} are algebraic structures that can distinguish framed links in certain lens spaces and framed \emph{virtual links} in thickened surfaces; see, for example, \cite{virtual}*{Sec.\ 3.2}. 

By Definition \ref{def:gl-rack}, GL-racks are precisely virtual racks in which all inner automorphisms $s_x$ are endomorphisms of virtual racks. Equivalently, a GL-rack $(X,s,\u)$ is a virtual rack in which the \emph{operator group} of $(X,s)$ identifies $x$ with $\u(x)$ for all $x\in X$; see \cite{fenn}*{Sec.\ 1.1}.
	\end{remark}
	
	The reader may have noticed that Definitions \ref{def:gl-rack} and \ref{def:gl-hom} are much simpler than the definitions originally given in \citelist{\cite{karmakar}*{Def.\ 3.1}\cite{bi}*{Def.\ 3.2}}. We will soon show that these definitions are equivalent. Before that, we consider several examples of GL-racks and state a short lemma.
	
	\begin{example}
		\label{ex:trivial}
		For all racks $R=(X,s)$, the identity map $\id_X$ is a GL-structure on $R$. Consequently, \emph{every} rack can be equipped with at least one GL-structure. 
	\end{example}
	
	\begin{example}[\cite{bi}*{Ex.\ 3.7}]\label{ex:constant}
		Given a permutation rack $P=(X,\sigma)\perm$, a GL-structure on $P$ is precisely a permutation $\u\in S_X$ such that $\u\sigma=\sigma\u$. Given such a $\u$, we say that $(P,\u)$ is a \emph{permutation GL-rack} or \emph{constant action GL-rack}, and we denote it by $(X,\sigma,\u)\perm$.
	\end{example}
	
	\begin{example}[\cite{bi}*{Ex.\ 3.6}]
		Let $G$ be a group, let $z\in Z(G)$ be a central element of $G$, and define $f:G\rightarrow G$ by $g\mapsto zg$. Then $(\Conj G, f)$ is a GL-quandle.
	\end{example}
	
	\begin{example}\label{ex:gl-dihedral}
		Let $n\geq 4$ be a multiple of $4$. Define four affine transformations ${f_{a,u}:\Z/n\Z\to\Z/n\Z}$ by $k\mapsto a+uk$ with $a\in\{0,n/2\}$ and $u\in \{1,1+n/2\}$. Each of these transformations is a GL-structure on the dihedral quandle $R_n$, and the translation defined by $k\mapsto k+1$ is a GL-rack isomorphism from $(R_n,f_{0,1+n/2})$ to $(R_n,f_{n/2,1+n/2})$. We later show that these are \emph{all} the possible GL-structures on $R_n$; see Proposition \ref{cor:dihedral}.
	\end{example}
	
	\begin{lemma}\label{lem:u-inv}
		If $\phi$ is a GL-rack homomorphism from $(R_1,\u_1)$ to $(R_2,\u_2)$, then $\phi\u_1\inv=\u_2\inv\phi$.
	\end{lemma}
	
	\begin{proof}
		By hypothesis, $\phi=\phi\u_1\u_1\inv=\u_2\phi\u_1\inv$. Applying $\u_2\inv$ on the left proves the claim. 
	\end{proof}
	
	\subsubsection{Bi-Legendrian racks}
Next, we reproduce the definition of GL-racks given in \citelist{\cite{karmakar}\cite{bi}}. Following Kimura \cite{bi}, we temporarily use the term \emph{bi-Legendrian racks} to distinguish them from GL-racks in the sense of Definition \ref{def:gl-rack}. We also define \emph{Legendrian racks}, which Ceniceros, Elhamdadi, and Nelson \cite{ceniceros} introduced in 2021.
	
	\begin{definition}\label{def:bi-leg}
		Given a rack $R=(X,s)$, a \emph{bi-Legendrian structure} on $R$ is a pair $(\u,\downcusp)$ of maps $\u,\downcusp:X\to X$ satisfying the following axioms for all elements $x\in X$:
		\begin{enumerate}[({L}1)]
			\item\label{L1} $\u\downcusp  s_x(x)=x=\downcusp\u  s_x(x)$.
			\item\label{L2} $\u s_x=s_x \u$,\quad and\quad$\downcusp s_x=s_x\downcusp$.
			\item\label{L3} $s_{\u(x)}=s_x=s_{\downcusp(x)}$.
		\end{enumerate}
		We call the triple $(R,\u,\downcusp)$ a \emph{bi-Legendrian rack}; if $\u=\downcusp$, then we also say that $(R,\u,\downcusp)$ is a \emph{Legendrian rack}. A \emph{bi-Legendrian rack homomorphism} from $(R_1,\u_1,\downcusp_1)$ to $(R_2,\u_2,\downcusp_2)$ is a rack homomorphism $\phi\in\Hom_\Rc(R_1,R_2)$ such that $\phi\u_1=\u_2\phi$ and $\phi \downcusp_1=\downcusp_2\phi$.
	\end{definition}

\begin{remark}\label{rmk:bi-quandle}
	Axiom \ref{L1} shows that, for all bi-Legendrian racks $(R,\u,\downcusp)$, the underlying rack $R$ is a quandle if and only if $\downcusp=\u\inv$.
\end{remark}

\subsection{Equivalence of definitions}
In this subsection, we show the equivalence of Definitions \ref{def:gl-rack} and \ref{def:bi-leg}. Specifically, we show that the category $\blr$ of bi-Legendrian racks is isomorphic to $\glr$.

\begin{prop}\label{prop:bi-leg-hom}
	If $(\u,\downcusp)$ is a bi-Legendrian structure on a rack $R=(X,s)$, then $\downcusp=\pi\inv\u\inv$, and $(R,\u)$ is a GL-rack. In particular, there is a forgetful functor $\Forg:\blr\to\glr$ sending $(R,\u,\downcusp)$ to $(R,\u)$. Moreover, $\Forg$ is fully faithful.
\end{prop}

\begin{proof}
	Axiom \ref{L1} states that \[\u\downcusp=\downcusp\u=\pi\inv,\] which is bijective, so $\u$ is bijective. It follows that $\downcusp=\pi\inv\u\inv$, as desired. Moreover, axioms \ref{L2} and \ref{L3} imply that $\u$ is a rack endomorphism that commutes with $s_x$ for all $x\in X$, so $(R,\u)$ is a GL-rack. Thus, we have a forgetful functor $\Forg:\blr\to\glr$, which is clearly faithful.
	
	Now, let $(R_1,\u_1,\downcusp_1)$ and $(R_2,\u_2,\downcusp_2)$ be bi-Legendrian racks. To show that $\Forg$ is full, we need to show that all GL-rack homomorphisms from $(R_1,\u_1)$ to $(R_2,\u_2)$ commute with the $\downcusp_i$'s. But this follows from the formula $\downcusp_i=\pi_{R_i}\inv\u_i\inv$, the inclusion $\Pi\inv\in Z(\Rc)$, and Lemma \ref{lem:u-inv}.
\end{proof}

\begin{prop}\label{prop:bi-is-gl}
	The forgetful functor $\Forg:\blr\to\glr$ is an isomorphism of categories. 
	Hence, the data and axioms of GL-racks are equivalent to those of bi-Legendrian racks. 
\end{prop}

\begin{proof}
	Since the functor $\Forg$ is fully faithful, it will suffice to show that $\Forg$ is bijective on objects. Injectivity follows from the first claim of Proposition \ref{prop:bi-leg-hom}.
	
	To show surjectivity, let $(R,\u)$ be a GL-rack with $R=(X,s)$, and define $\downcusp:=\pi\inv\u\inv$ as in Proposition \ref{prop:bi-leg-hom}. 
	It will suffice to show that $(R,\u,\downcusp)$ is a bi-Legendrian rack. To that end, fix $x\in X$. Since $\Pi\inv\in Z(\Rc)$,
	\[
	\ud s_x(x)=\u\pi\inv\u\inv \pi(x)=\u\u\inv\pi\inv\pi(x)=x
	\]
	and, similarly, $\du s_x(x)=x$. This verifies axiom \ref{L1}. Lemma \ref{lem:u-inv} and the inclusion $\Pi\inv\in Z(\Rc)$ imply that $\downcusp s_x=s_x\downcusp$, which verifies axiom \ref{L2}. Since $\u$ is a GL-structure, \[s_{\u(x)}\u=\u s_x=s_x\u.\]
	From the bijectivity of $\u$, we obtain $s_{\u(x)}=s_x$. 
	By construction, $\downcusp\in \Aut R$, so we similarly obtain $s_{\downcusp(x)}=s_x$. 
	This verifies axiom \ref{L3}. Hence, $(R,\u,\downcusp)$ is a bi-Legendrian rack, as desired.
\end{proof}
	
	In light of these results, we will henceforth call bi-Legendrian racks GL-racks except when specifically citing the axioms in Definition \ref{def:bi-leg} or denoting GL-racks as quadruples with $\downcusp$. The previous two propositions make computing and classifying GL-racks significantly easier; see Section \ref{subsec:classn}. As a bonus, they also yield the following converse of \cite{ceniceros}*{Rem.\ 2}.
	
	\begin{cor}\label{cor:leg-quandles}
		A GL-rack $(R,\u)$ is a Legendrian rack if and only if $\pi_R=\u^{-2}$. In this case, $R$ is a quandle if and only if $\u$ is an involution.
	\end{cor}
	
	\subsection{GL-racks as an algebraic theory}\label{sec:ftrs}
	In this subsection, we adapt the discussion in Section \ref{sec:rack-thy} to GL-racks, viewed as an algebraic theory with two binary operations $s_-(-)$ and $s\inv_-(-)$ and two unary operations $\u$ and $\u\inv$. 
	
	\subsubsection{Examples of quotient GL-racks}
	First, we adapt the definitions of Section \ref{sec:rack-thy} from $\Rc$ to $\glr$.
	
	\begin{definition}
		Given a GL-rack $R=(X,s,\u)$, the \emph{associated GL-quandle} of $R$ is the GL-quandle $R\qnd:=R/{\sim}$, where $\sim$ is the smallest congruence on $R$ such that $s_x(x)\sim x$ for all $x\in X$. Similarly, the \emph{medialization} or \emph{abelianization} of $R$ is the medial GL-rack $R\med:=R/{\sim}$, where $\sim$ is the smallest congruence on $R$ such that equation (\ref{eq:medial-2}) holds for all $x,y,z,a\in X$. 
	\end{definition}
	
	\subsubsection{Free GL-racks}
	As with racks, the \emph{free GL-rack} $\Free_\glr(X)$ on a set $X$ exists and is uniquely (up to isomorphism) characterized by a universal property; see \cite{karmakar}*{Prop.\ 4.2}. Namely, there exists a set map $i:X\to\Free_\glr(X)$ such that for all GL-racks $R=(Y,t,\u)$ and all set maps $\phi_0:X\to Y$, there exists a unique GL-rack homomorphism $\phi\in\Hom_\glr(\Free_\glr(X),R)$ such that $\phi i=\phi_0$. That is, the following diagram commutes:
	\begin{equation}\label{diagram:free-glr}\begin{tikzcd}
			X & {\operatorname{Free}_{\mathsf{GLR}}(X)} \\
			& Y
			\arrow["i", from=1-1, to=1-2]
			\arrow["{\phi_0}"', from=1-1, to=2-2]
			\arrow["\phi", dashed, from=1-2, to=2-2]
	\end{tikzcd}\end{equation}
	Karmakar, Saraf, and Singh \cite{karmakar}*{Sec.\ 4} gave an explicit set-theoretic construction of $\Free_\glr(X)$. 
	Here, we simplify their construction using Definition \ref{def:gl-rack} and and Proposition \ref{prop:bi-is-gl}.
	\begin{definition}
		Let $X$ be a set. We define the \emph{free GL-rack on $X$} as follows. If $X=\emptyset$, let $\Free_\glr(X)$ be the empty GL-rack.
		Else, let the \emph{universe of words generated by $X$} be the set $W(X)$ such that $X\subset W(X)$ and $s_y(x),\,s_y\inv(x), \,\u(x),\,\u\inv(x)\in W(X)$ for all $x,y \in W(X)$.
		Let $V(X)$ be the set of equivalence classes of elements of $W(X)$ modulo the congruence generated by the following relations for all $x,y,z \in W(X)$:
		\begin{itemize}
			\item $s_y\inv s_y(x)  \sim  s_ys_y\inv(x)\sim x\sim \u\u\inv(x)\sim \u\inv\u(x)$.
			\item $s_zs_y(x)\sim s_{s_z(y)}s_z(x)$.
			\item $s_{\u(y)}\u(x) \sim\u s_y(x) \sim s_y\u(x)$.
		\end{itemize}
		Thus, we have a rack structure $s:V(X)\to S_{V(X)}$ on $V(X)$ and a GL-structure $\u\in S_{V(X)}$ on $(V(X),s)$. So, we define $\Free_\glr(X)$ to be the GL-rack $(V(X), s, \u)$.
	\end{definition}
	
	\begin{remark}\label{remark:free-bi}
		By way of Proposition \ref{prop:bi-leg-hom}, for all $x\in V(X)$, we can consider $\downcusp(x):=\pi\inv\u\inv(x)$ as an element of $V(X)$.
	\end{remark}
	
	\subsubsection{Free objects in subtheories of GL-racks}

		Since $\glr$ is the category of models of the algebraic theory of GL-racks in $\Set$, the free GL-rack $L=\Free_\glr(X)$ on the one-element set $X=\{x\}$ is a \emph{strong generator} or \emph{separator} for $\glr$; see, for example, \cite{alg-theories}*{Prop.\ 3.3.3}. 

Since GL-quandles, Legendrian racks, and Legendrian quandles are subtheories of the algebraic theory of GL-racks, we can similarly consider free objects that strongly generate the corresponding subcategories of models by taking quotients of $L$ by congruences. 

\begin{definition}
	Given a set $Y$, let $R$ be the free GL-rack on $Y$. We define the \emph{free GL-quandle} on $Y$ to be $R\qnd$. 
	
	Using Corollary \ref{cor:leg-quandles}, we similarly define the \emph{free Legendrian rack} on $Y$ to be the quotient of $R$ by the smallest congruence $\sim$ such that $s_y(y)\sim\u^{-2}(y)$ for all $y\in V(Y)$, and we define the \emph{free Legendrian quandle} on $Y$ to be $(R/{\sim})\qnd$. 
\end{definition}

By the above discussion, we have the following.

	\begin{prop}\label{rmk:generators}
		Let $L$ be the free GL-rack on one element. Then $L\qnd$, $L/{\sim}$, and $(L/{\sim})\qnd$ are strong generators of $\glq$, the category of Legendrian racks, and the category of Legendrian quandles, respectively. 
Moreover, as sets,
\[L\qnd=L/{\sim}=\{\u^k (x)\mid k\in\Z\}, \quad \text{and } (L/{\sim})\qnd=\{x,\u (x)\}.\]
	\end{prop}
	
	\subsubsection{The free GL-rack on one generator}
	We prove an analogue of Example \ref{ex:free-one} for GL-racks.
	
	\begin{prop}\label{prop:free-gl-one}
		The free GL-rack on one generator is canonically isomorphic to the permutation GL-rack $L:=(\Z^2,\sigma,\u_0)\perm$, where $\sigma(m,n)=(m+1,n)$ and $\u_0(m,n)=(m,n+1)$ for all $m,n\in\Z$.
	\end{prop}
	
	\begin{proof}
	Let $X=\{x\}$. We will show that $L$ satisfies the universal property of $\Free_\glr(X)$ with $i:X\to\Z^2$ defined by $x\mapsto(0,0)$. To that end, let $(R,\u)$ be a GL-rack with $R=(Y,t)$, and let $\phi_0:X\to Y$ be a set map. 
	Define $\phi:\Z^2\to Y$ by \[(m,n)\mapsto \u^n \pi_R^m \phi_0(x).\]  
	
	First, we show that $\phi$ is a GL-rack homomorphism from $L$ to $(R,\u)$. To that end, fix elements $(k,\ell),(m,n)\in\Z^2$. Denote the underlying rack structure on $L$ by $s:\Z^2\to S_{\Z^2}$ with $s_{(m,n)}:=\sigma$ for all $(m,n)\in\Z^2$. By part \ref{ii} of Proposition \ref{prop:rack-center}, 
	\[
	\phi s_{(k,\ell)}(m,n)=\phi\sigma(m,n)=\phi(m+1,n)=\u^n\pi_R^{m+1}\phi_0(x)=\u^n t^{m+1}_{\phi_0(x)}\phi_0(x).
	\]
	On the other hand, we apply bi-Legendrian rack axiom \ref{L3}, part \ref{iii} of Proposition \ref{lem:theta-inv}, and part \ref{ii} of Proposition \ref{prop:rack-center} to compute
	\[
	t_{\phi(k,\ell)}\phi(m,n)=t_{\u^\ell \pi_R^k\phi_0(x)}\u^n \pi_R^m\phi_0(x)= t_{\phi_0(x)}\u^n t^m_{\phi_0(x)}\phi_0(x)=\u^n t^{m+1}_{\phi_0(x)}\phi_0(x)=\phi s_{(k,\ell)}(m,n),
	\]
	so $\phi$ is a rack homomorphism. Moreover, \[\phi\u_0(m,n)=\phi(m,n+1)=\u^{n+1}\pi_R^m\phi_0(x)=\u\phi(m,n),\] so $\phi$ is a GL-rack homomorphism.
	
	Next, we show that diagram (\ref{diagram:free-glr}) commutes and that $\phi$ is unique. Indeed, we have commutativity because \[\phi i(x)=\phi(0,0)=\u^0\pi^0_R\phi_0(x)=\phi_0(x).\] Now, let $\psi\in\Hom_\glr(L,(R,\u))$ be any GL-rack homomorphism such that $\psi i(x)=\phi_0(x)$, and fix $(m,n)\in\Z^2$. By part \ref{ii} of Proposition \ref{prop:rack-center},
	\[
	\u^{-n} \pi_R^{-m} \psi(m,n)= \u^{-n} t^{-m}_{\psi(m,n)}\psi(m,n) =\psi\u_0^{-n}s^{-m}_{(m,n)}(m,n)=\psi(0,0)=\psi i(x)=\phi_0(x)
	\]
	and, similarly, \[\u^{-n} \pi_R^{-m} \phi(m,n)=\phi_0(x).\] Since $\u$ and $\pi_R$ are bijections, we obtain $\psi(m,n)=\phi(m,n)$. Since the element $(m,n)\in\Z^2$ was arbitrary, this shows that $\phi$ is unique. Hence, $L$ satisfies the universal property of $\Free_\glr(X)$.
\end{proof}

	A similar proof yields the following generalization of Proposition \ref{prop:free-gl-one}; the author is grateful to the anonymous referees for observing this.
	
	\begin{thm}\label{thm:free-arbitrary}
		Let $X$ be a set, and let $\Free_{\Rc}(X)=(F(X),t)$ be the free rack generated by $X$. Then the free GL-rack generated by $X$ can be constructed as the product of $\Free_{\Rc}(X)$ and the trivial quandle $(\Z,\id_\Z)\perm$: 
		\[
		\Free_\glr(X)\cong (F(X)\times\Z,s,\u),\quad s_{(y,n)}(x,m):=(t_y(x),m),\quad \u(x,m):=(x,m+1).
		\]
	\end{thm}

	\subsubsection{Mediality and Hom GL-racks}
	
	Several quandle-theoretic invariants of smooth links can be enhanced using Corollary \ref{cor:hom-racks}; see, for example, \citelist{\cite{Hom}\cite{enhancements}}. 
	This motivates the following analogue of Corollary \ref{cor:hom-racks} for GL-racks; for applications, see Example \ref{ex:bihoms}, Theorem \ref{thm:scmc}, and Corollary \ref{prop:medial-enhance}.
	
	\begin{thm}\label{thm:homsets}
		Let $R_1:=(X,s,\u_1)$ and $R_2:=(Y,t,\u_2)$ be GL-racks, and suppose that $R_2$ is medial. Then $H:=\Hom_\glr(R_1,R_2)$ is a subrack of $\Hom_\Rc((X,s),(Y,t))$ equipped with its medial rack structure $\Tilde{t}$ from Corollary \ref{cor:hom-racks}, and $(H,\Tilde{t}|_H)$ has a canonical GL-structure $\u:H\rightarrow H$ defined by $f\mapsto \u_2 f$. 
		In particular, if $R_1$ or $R_2$ is a GL-quandle, then so is $(H,\Tilde{t}|_H,\u)$.
	\end{thm}
	
	\begin{proof}
		By Proposition \ref{prop:aut-inn-dual} and Lemma \ref{lem:u-inv}, the binary operations $s_-(-)$ and $s\inv_-(-)$ and the unary operations $\u$ and $\u\inv$ are all homomorphisms in the algebraic theory of medial GL-racks. Therefore, medial GL-racks are a commutative algebraic theory, so the claim follows directly from Proposition \ref{prop:hom-alg-thy} and Corollary \ref{cor:hom-racks}.
	\end{proof}
	
	In particular, Theorem \ref{thm:homsets} yields the following theoretical enhancement of the GL-rack coloring numbers studied in \citelist{\cite{ceniceros}\cite{bi}\cite{karmakar}}. This medial GL-rack-valued invariant is inspired by similar enhancements of quandle colorings of smooth links in \citelist{\cite{Hom}\cite{enhancements}}.
	
	\begin{cor}\label{prop:medial-enhance}
		Let $\GG(\Lambda)$ be the fundamental GL-rack of a Legendrian link $\Lambda$ (see \cite{karmakar}). Then for all medial GL-racks $M$, the isomorphism type of $\Hom_\glr(\GG(\Lambda),M)$ as a medial GL-rack is an invariant of $\Lambda$.
	\end{cor}
	
	\section{Group-theoretic aspects of GL-racks}\label{subsec:classn}
	In a previous version of \cite{karmakar}*{Sec.\ 3}, Karmakar, Saraf, and Singh posed the following question: what are all the possible GL-structures on a given rack? In this section, we answer this question and classify various infinite families of GL-racks. As applications, we discuss automorphism groups of GL-racks and compute these groups for all GL-racks whose underlying racks are dihedral quandles. 
	
	\subsection{Characterizing GL-structures}
	Given a rack $R=(X,s)$, let $U_R$ be the set of GL-structures $\u:X\to X$ on $R$, and define an equivalence relation $\sim$ on $U_R$ by identifying $\u_1\sim\u_2$ if and only if $(R,\u_1)\cong (R,\u_2)$. 
	The following theorem completely characterizes $U_R$ and $U_R/{\sim}$. 
	Given a group $G$ and a subset or element $H$ of $G$, we will denote the centralizer of $H$ in $G$ by $C_G(H)$.
	
	\begin{thm}\label{cor:number}
		Given a rack $R=(X,s)$, define $C:=C_{\Aut R}(\operatorname{Inn}R)$. Then $U_R=C$, and $U_R$ is a normal subgroup of $\Aut R$. Furthermore, $\u_1\sim\u_2$ if and only if $\u_1$ and $\u_2$ are conjugate in $\Aut R$. 
		In particular, if $\Aut R$ is abelian, then $U_R =U_R/{\sim}=\Aut R$.
	\end{thm}
	
	\begin{proof}
		The claim that $U_R=C$ is a restatement of Definition \ref{def:gl-rack}. It is straightforward to verify that $\Inn R$ is a normal subgroup of $\Aut R$, so $C$ is normal in $\Aut R$. 

On the other hand, given two GL-structures $\u_1,\u_2\in U_R=C$, a map $\phi:X\to X$ is a GL-rack isomorphism from $(R,\u_1)$ to $(R,\u_2)$ if and only if $\phi\in\Aut R$ and $\phi\u_1=\u_2\phi$. In other words, $(R,\u_1)\cong (R,\u_2)$ if and only if $\u_1$ and $\u_2$ are conjugate in $\Aut R$, as claimed.

Finally, suppose that $\Aut R$ is abelian. Then $\Aut R=C=U_R$, and each element of $\Aut R$ constitutes its own conjugacy class. Hence, $U_R=U_R/{\sim}$.
	\end{proof}
	
	One may ask whether $U_R$ is conjugacy-closed in $\Aut R$, that is, whether $\u_1\sim\u_2$ in $U_R/{\sim}$ if and only if $\u_1$ and $\u_2$ are conjugate in $U_R$. We will give a negative answer later; see Remark \ref{rmk:not-conj-closed}.
	
	\begin{cor}
		For all racks $R$, we have $U_R=U_{R^{\operatorname{op}}}$ and $U_R/{\sim}=U_{R^{\operatorname{op}}}/{\sim}$.
	\end{cor}
	
	\begin{proof}
		This follows immediately from Proposition \ref{prop:aut-inn-dual} and Theorem \ref{cor:number}.
	\end{proof}
	
	\subsection{Classification of GL-racks}
	In this subsection, we use Theorem \ref{cor:number} to classify certain families of GL-racks.
	
		\subsubsection{Computer classification}
	
	As an application of Theorem \ref{cor:number}, we used \texttt{GAP} \cite{GAP4} to classify all GL-racks up to order $8$ up to isomorphism. Our \texttt{GAP} program and raw data are available in a GitHub repository \cite{my-code}. Our program builds upon the work of Vojtěchovský and Yang \cite{library} in 2019, who classified racks up to order $11$ \cite{library-site}. Table \ref{tab:tab1} enumerates our data.
	
	\begin{table}[h]
		\caption{Enumeration of GL-racks up to order $8$.}
		\label{tab:tab1}
		\centering
		\begin{tabular}{l|lllllllll}
			Order      & $0$ & $1$ & $2$ & $3$  & $4$  & $5$   & $6$    & $7$  & $8$   \\ \hline
			GL-racks  & $1$ & $1$ & $4$ & $13$ & $62$ & $308$ & $2132$ & $17268$ & $189373$ \\
			Medial GL-racks & $1$ & $1$ & $4$ & $13$ & $61$ & $298$ & $2087$ & $16941$ & $187160$ \\
			GL-quandles & $1$ & $1$ & $2$ & $6$  & $19$ & $74$  & $353$  & $2080$ & $16023$ \\
			Medial GL-quandles & $1$ & $1$ & $2$ & $6$  & $18$ & $68$  & $329$  & $1965$ & $15455$  
		\end{tabular}
	\end{table} 
	
	\subsubsection{Permutation GL-racks} First, we classify GL-structures on permutation racks.
	
	\begin{prop}\label{cor:perms}
		Let $X$ be a set, let $\sigma\in S_X$, and let $P$ be the permutation rack $(X,\sigma)\perm$. Then $U_P=C_{S_X}(\sigma)=\Aut P$, and $U_P/{\sim}$ is the set of conjugacy classes of $C_{S_X}(\sigma)$.
	\end{prop}
	
	\begin{proof}
		An automorphism of $P$ is precisely a permutation $\phi\in S_X$ such that $\phi\sigma=\sigma\phi$. Therefore, $\Aut P=C_{S_X}(\sigma)$. On the other hand, $\operatorname{Inn}P=\langle\sigma\rangle$, so Theorem \ref{cor:number} states that \[U_P=C_{C_{S_X}(\sigma)}(\langle\sigma\rangle)=C_{C_{S_X}(\sigma)}(\sigma)=C_{S_X}(\sigma)=\Aut P,\]as desired. Combined with Theorem \ref{cor:number}, these equalities imply the last part of the claim.
	\end{proof}
	
	\begin{example}\label{ex:perms}
		For all trivial quandles $P=(X,\id_X)\perm$, we have $C_{S_X}(\id_X)=S_X$. So, Proposition \ref{cor:perms} states that $U_P=S_X$, and $U_P/{\sim}$ is the set of conjugacy classes of $S_X$.
	\end{example}
	
	\begin{example}
		Let $X:=\{1,2,\dots,n\}$, let $\sigma\in S_n$ be an $n$-cycle, and let $P:=(X,\sigma)\perm$. Then  $C_{S_n}(\sigma)$ is the cyclic subgroup $\langle \sigma\rangle\cong \Z/n\Z$ of $S_n$; see, for example, \cite{dummit}*{p.\ 127}. Since $\langle \sigma\rangle$ is abelian, Proposition \ref{cor:perms} implies that $U_P/{\sim}=U_P=\langle \sigma\rangle$.
	\end{example}
	
	\begin{example}\label{ex:free-one-aut}
		Let $F$ be the free rack on one element. For all $n\in\Z$, let $\tau_n:\Z\to\Z$ be the translation defined by $k\mapsto k+n$. Recall from Example \ref{ex:free-one} that $F$ is isomorphic to the permutation rack $(\Z,\tau_1)\perm$. 
		By Proposition \ref{cor:perms}, \[U_F=\Aut F=C_{S_\Z}(\tau_1)=\{\tau_n\mid n\in\Z\}\cong \Z.\] In particular, $\Aut F$ is abelian, so $U_F=U_F/{\sim}$. 
		In other words, there are infinitely many GL-structures on $F$, all of which are translations of $\Z$ and none of which yield isomorphic GL-racks.
	\end{example}
	
	\subsubsection{Conjugation GL-quandles}
	Next, we use Theorem \ref{cor:number} to classify GL-structures on conjugation quandles of centerless or abelian groups. Given a group $G$, let $\Aut_\Grp G$ and $\operatorname{Inn}_\Grp G$ denote the automorphism group and inner automorphism group of $G$, respectively. 
	
	\begin{prop}\label{cor:conj}
		Let $G$ be a group, and let $Q:=\operatorname{Conj}G$. If $G$ is abelian, then $U_Q=S_G$, and $U_Q/{\sim}$ is the set of conjugacy classes of $S_G$. On the other hand, if $G$ is centerless, then $U_Q=\{\id_G\}$.
	\end{prop}
	
	\begin{proof}
		If $G$ is abelian, then $Q$ is a trivial quandle, so the first claim follows from Example \ref{ex:perms}. 
		Next, recall a result of Elhamdadi, Macquarrie, and Restrepo \cite{dihedral}*{Thm.\ 2.3} that $\operatorname{Inn}Q= \operatorname{Inn}_\Grp G$ for all groups $G$. 
		Also, recall a result of Bardakov, Nasybullov, and Singh \cite{centerless}*{Cor.\ 2} that $\Aut Q=\Aut_\Grp G$ if and only if $G$ is centerless. In this case, Theorem \ref{cor:number} states that \[U_Q= C_{\Aut_\Grp G}(\operatorname{Inn}_\Grp G),\] which is the group of central automorphisms of $G$. However, this group is trivial when $G$ is centerless; see, for example, \cite{robinson}*{p.\ 410}.
	\end{proof}
	
	\begin{example}
		Let $n\geq 3$ be an integer or $n=\infty$. Then the symmetric group $S_n$ is centerless, so Proposition \ref{cor:conj} states that the only GL-structure on $\operatorname{Conj}S_n$ is $\id_{S_n}$. Similarly, if $n\geq 4$, then the alternating group $A_n$ is centerless, so the only GL-structure on $\operatorname{Conj}A_n$ is $\id_{A_n}$.
	\end{example}
	
	\subsubsection{Takasaki GL-kei}
	Next, we classify GL-structures on a certain family of Takasaki kei and, as a consequence, all dihedral quandles of odd order. 
	
	We first recall a classification result of Bardakov, Dey, and Singh \cite{takasaki-aut}*{Thm.\ 4.2}. For all abelian additive groups $A$ without $2$-torsion, $\Aut T(A)$ is isomorphic to the holomorph \[G:=A\rtimes\Aut_\Grp A\] of $A$. 
Under this identification, 
$\operatorname{Inn}T(A)$ is the semidirect product \[H:=2A\rtimes\{\pm \id_A\}\leq G,\] where $-\id_A$ denotes inversion. We prove the following result using these identifications.
	
	\begin{prop}\label{cor:takasaki}
		If $A$ is an abelian additive group without $2$-torsion, then the only GL-structure on the Takasaki kei $T(A)$ is $\id_A$.
	\end{prop}
	
	\begin{proof}
		First, note that for all automorphisms $\psi\in \Aut_\Grp A$ such that $\psi(2a)=2a$ for all $a\in A$, we have $2(\psi(a)-a)=0$. Since $A$ is $2$-torsion-free, it follows that $\psi(a)=a$, so $\psi=\id_A$. Therefore, by Theorem \ref{cor:number}, 
		it will suffice to show that 
		\[
		C_G(H)\subseteq\{(0,\psi)\in G: \psi|_{2A}=\id_{2A}\}
		\]
		since, as we just observed, the right-hand side is the trivial subgroup of $G$. 
		To that end, a direct computation shows that conjugation in $G$ is given by
		\begin{equation*}
			(a,\psi)(b,\phi)(a,\psi)\inv=(a+\psi(b)-\psi\phi\psi\inv(a),\psi\phi\psi\inv).
		\end{equation*}
		For all $(b,\phi)\in H$, we have $\phi=\pm\id_A$. It follows that, for all $(a,\psi)\in C_G(H)$ and $(b,\phi)\in H$, 
		\[
		(b,\phi)=(a,\psi)(b,\pm\id_A)(a,\psi)\inv=(a+\psi(b)\mp a,\pm\id_A).
		\]
		Taking $\phi:=+\id_A$ yields $\psi(b)=b$; since this equality holds for all $(b,\phi)\in H$ (and, hence, for all $b\in 2A$), we obtain ${\psi|_{2A}=\id_{2A}}$, as desired. Therefore, taking $\phi:=-\id_A$ yields $b=2a+b$, so $2a=0$. Since $A$ has no $2$-torsion, it follows that $a=0$, as desired. 
	\end{proof}
	
	\begin{example}\label{ex:dihedral-odd}
		For all odd integers $n\geq 3$, the cyclic group $\Z/n\Z$ is $2$-torsion-free. By Proposition \ref{cor:takasaki}, the only GL-structure on the dihedral quandle $R_n$ of order $n$ is $\id_{\Z/n\Z}$.
	\end{example}
	
	\subsubsection{Dihedral GL-quandles}
	Without the assumption that $A$ is $2$-torsion-free, there are infinitely many counterexamples to Proposition \ref{cor:takasaki}. 
To show this, we will complete the classification of GL-structures on dihedral quandles $R_n$ using the following results of Elhamdadi, Macquarrie, and Restrepo \cite{dihedral}*{Thms.\ 2.1 and 2.2}. The automorphism group $\Aut R_n$ is the group of affine transfomations of $\Z/n\Z$. Thus, $\Aut R_n$ is isomorphic to the holomorph \[G:=\Z/n\Z\rtimes(\Z/n\Z)^\times\]of $\Z/n\Z$. 
Under this identification, 
$\operatorname{Inn}R_n$ is the semidirect product \[H:=2\Z/n\Z\rtimes\{\pm 1\}\leq G,\] which is isomorphic to the dihedral group $D_{n/2}$ of order $n$.
	
	The following result strengthens Example \ref{ex:gl-dihedral}, and we state it in terms of the above identifications. Our proof uses the fact that conjugation in $G$ is given by
	\begin{equation}\label{eq:dihedral}
		(a,u)(b,v)(a,u)\inv=(a+ub,uv)(-u\inv a,u\inv)=(ub+(1-v)a,v).
	\end{equation}
	
	\begin{prop}\label{cor:dihedral}
		For all even integers $n\geq 2$, the GL-structures on the dihedral quandle $R_n$ are
		\begin{equation}\label{eq:urn}
			U_{R_n}=\begin{cases}
				\{0,n/2\}\rtimes\{1\}& \text{if }4\nmid n,\\
				\{0,n/2\} \rtimes \{1,1+n/2\} & \text{if }4\mid n.
			\end{cases}
		\end{equation}
		If $4\nmid n$, then $U_{R_n}=U_{R_n}/{\sim}$, so $|U_{R_n}/{\sim}|=2$. 
		Otherwise, the only elements of $U_{R_n}$ that are identified in $U_{R_n}/{\sim}$ are $(0,1+n/2)$ and $(n/2,1+n/2)$, so $|U_{R_n}/{\sim}|=3$.
	\end{prop}
	
	\begin{proof}
		Theorem \ref{cor:number} states that, to prove equation (\ref{eq:urn}), it will suffice to show that the right-hand side equals $C_G(H)$. For all elements $(a,u)\in C_G(H)$ and $(b,v)\in H$, the right-hand side of equation (\ref{eq:dihedral}) equals $(b,v)$. In particular, \[b=ub+(1-v)a\] for all $b\in 2\Z/n\Z$, so taking $v:=1$ yields $u\in\{1,1+n/2\}$. However, $1+n/2\in(\Z/n\Z)^\times$ if and only if $4\mid n$, as desired. Since $ub=b$, taking $v:=-1$ yields $2a=0$. Hence, $a\in\{0,n/2\}$, as desired. 
		This shows that $C_G(H)$ is a subset of the right-hand side of equation (\ref{eq:urn}), and verifying the opposite containment is straightforward.
		
		We now prove the second claim. Since $(0,1)$ is the identity element of $G$, it is not conjugate to any other element of $U_{R_n}$. If $4\nmid n$, then we are done. 
		
		Otherwise, let $(b,v),(c,w)\in U_{R_n}$. Then $(b,v)$ and $(c,w)$ are conjugate in $G$ if and only if there exists an element $(a,u)\in G$ such that $(c,w)$ equals the right-hand side of equation (\ref{eq:dihedral}). In particular, $w=v$. It follows that neither $(0,1+n/2)$ nor $(n/2,1+n/2)$ is conjugate to $(n/2,1)$ in $G$. 
		On the other hand, taking $(b,v):=(0,1+n/2)$ and $(a,u):=(1,1)$ in equation (\ref{eq:dihedral}) shows that $(0,1+n/2)$ and $(n/2,1+n/2)$ are conjugate in $G$, so the proof is complete.
	\end{proof}
	
	\begin{remark}\label{rmk:not-conj-closed}
		If $n\geq 4$ is a multiple of $4$, then equation (\ref{eq:urn}) shows that $U_{R_n}\cong\Z/2\Z\times\Z/2\Z$, so $U_{R_n}$ is abelian.
		It follows that $(0,1+n/2)$ and $(n/2,1+n/2)$ are conjugate in $G$ but not in $U_{R_n}$. 
		Hence, the condition in Theorem \ref{cor:number} that $\u_1\sim\u_2$ in $U_R/{\sim}$ does not imply conjugacy in $U_R$.
	\end{remark}
	
	\subsection{Automorphism groups of GL-racks}
Given a GL-rack $(R,\u)$, let $\Aut_\glr (R,\u)$ denote its group of GL-rack automorphisms. The following characterization is simply a restatement of the definition of GL-rack automorphisms.

\begin{prop}\label{thm:aut-grp}
	For all GL-racks $(R,\u)$, we have $\Aut_\glr (R,\u)=C_{\Aut R}(\u)$.
\end{prop}

\begin{example}\label{ex:auts-of-free-gl-one}
	Let $L$ be the free GL-rack on one element, and identify $L=(\Z^2,\sigma,\u_0)\perm$ as in Proposition \ref{prop:free-gl-one}. 
	It is straightforward to show that GL-rack endomorphisms of $L$ are precisely translations of the form $(m,n)\mapsto(m+k,n+\ell)$ for some $(k,\ell)\in\Z^2$. Since all maps of this form are permutations of $\Z^2$, the mapping $(k,\ell)\mapsto \u_0^\ell \sigma^k$ is a group isomorphism from $\Z^2$ to $\Aut_\glr(L)$.
\end{example}

\subsubsection{Automorphisms of dihedral GL-quandles} We classify automorphism groups of GL-racks whose underlying racks are dihedral quandles $R_n$ of order $n$. Once again, we use a result of Elhamdadi, Macquarrie, and Restrepo \cite{dihedral}*{Thm.\ 2.1} to identify $\Aut R_n\cong \Z/n\Z\rtimes (\Z/n\Z)^\times$. 
	
	\begin{prop}
		Let $n\geq 2$ be an integer, let $R_n$ be the dihedral quandle of order $n$, and let $f_{b,v}:\Z/n\Z\to\Z/n\Z$ defined by $k\mapsto b+vk$ be a GL-structure on $R_n$; see Proposition \ref{cor:dihedral}. Let $G:=\Aut_\glr (R_n,f_{b,v})$. Then
\[
G\cong \begin{cases}
	2\Z/n\Z\rtimes(\Z/n\Z)^\times& \text{if }4\mid n\text{ and }v=1+\frac{n}{2},\\
	\Z/n\Z\rtimes(\Z/n\Z)^\times& \text{otherwise}.
\end{cases}
\]
In the latter case, $G=\Aut R_n$.
	\end{prop}
	
	\begin{proof}
		By Proposition \ref{thm:aut-grp}, $G=C_{\Aut R_n}((b,v))$, so for all $(a,u)\in \Aut R_n$, we have $(a,u)\in G$ if and only if $(b,v)$ equals the right-hand side of equation (\ref{eq:dihedral}). 
		Certainly, if $(b,v)=(0,1)$, then 
		$G=\Aut R_n$, as claimed. 
		Otherwise, $n$ is even by Example \ref{ex:dihedral-odd}.
		By Proposition \ref{cor:dihedral}, it suffices to only consider the cases that $(b,v)=(n/2,1)$ and $(b,v)=(0,1+n/2)$.
		
		If $(b,v)=(n/2,1)$, then $(a,u)\in G$ if and only if $0=(u-1)(n/2)$ in $\Z/n\Z$. 
		But $u-1$ is even for all elements $u\in (\Z/n\Z)^\times$, so this equation always holds. In other words, 
		\emph{every} element $(a,u)\in \Aut R_n$ centralizes $(b,v)$, so $G=\Aut R_n$.
		
		If ${(b,v)=(0,1+n/2)}$, then $(a,u)\in G$ if and only if $0=(-n/2)a$ in $\Z/n\Z$. This is true if and only if $a$ is even, and there are no restrictions placed on $u$. It follows that $G=2\Z/n\Z\rtimes(\Z/n\Z)^\times$, which completes the proof.
	\end{proof}
	
	\section{Categorical aspects of GL-racks}\label{sec:cat-equiv}
	In this section, we compute the centers of $\glr$ and several of its full subcategories. We also prove a surprising equivalence of categories between $\Rc$ and $\glq$.
	\subsection{Categorical centers}
	 Recall that the \emph{center} of a category $\mathcal{C}$ is the commutative monoid $Z(\mathcal{C})$ of natural endomorphisms of the identity functor $\mathbf{1}_\mathcal{C}$. In 2018, Szymik \cite{center}*{Thms.\ 5.4 and 5.5} computed that $Z(\Rc)=\langle\Pi\rangle\cong \Z$ and $Z(\Qn)\cong \{1\}$. 
In this section, we similarly compute the centers of $\glr$, $\glq$, and the categories of Legendrian racks and Legendrian quandles.

\begin{thm}
	Let $\Pi$ be the collection of canonical automorphisms $\pi_R$ of racks $R$, and let $\u$ be the collection of all GL-structures on racks. Then we have the following:
	\begin{enumerate}[({Z}1)]
		\item\label{Z1} The center $Z(\glr)$ is the free abelian group $\langle \Pi,\u\rangle\cong\Z^2$ generated by $\Pi$ and $\u$.
		\item\label{Z2} The centers of $\glq$ and the category of Legendrian racks are each the free group $\langle\u\rangle\cong\Z$.
		\item\label{Z3} The center of the category of Legendrian quandles is the group $\langle \u\mid\u^{2}=1\rangle\cong\Z/2\Z$.
	\end{enumerate}
\end{thm}

\begin{proof}
	To prove claim \ref{Z1}, let $\eta$ be a natural endomorphism of $\mathbf{1}_\glr$. By definition, $\eta$ is contained in the center $Z(\glr)$ if and only if, for all GL-racks $R_1=(X,s,\u_1)$ and $R_2=(Y,t,\u_2)$ and GL-rack homomorphisms $\phi\in\Hom_\glr(R_1,R_2)$, the following diagram commutes:
	\[\begin{tikzcd}
		X & X \\
		Y & Y
		\arrow["{\eta_{R_1}}", from=1-1, to=1-2]
		\arrow["\phi"', from=1-1, to=2-1]
		\arrow["\phi", from=1-2, to=2-2]
		\arrow["{\eta_{R_2}}"', from=2-1, to=2-2]
	\end{tikzcd}\]
	To see that $\langle\Pi,\u \rangle\subseteq Z(\glr)$, note by the definition of a GL-rack homomorphism that taking $\eta:=\u^n$ with $n\in\Z$ makes the diagram commute. If we take $\eta:=\Pi^m$ with $m\in\Z$, then the diagram commutes because $\Pi$ generates $Z(\Rc)$. In particular, taking $R_2:=R_1$ and $\phi:=\u_1^n$ shows that $\Pi^m$ and $\u^n$ commute in $Z(\glr)$. 
	
	To see that $Z(\glr)\subseteq\langle\Pi,\u \rangle\cong\Z^2$, take $R_1$ to be the free GL-rack on one element $(0,0)$ with its identification $R_1=(\Z^2,\sigma,\u_0)\perm$ from Proposition \ref{prop:free-gl-one}. 
	Let $\phi_0:\{(0,0)\}\to Y$ be a set map, let $\phi:\Z^2\to Y$ be the induced GL-rack homomorphism from the universal property of free GL-racks, and fix $\eta\in Z(\glr)$. By Example \ref{ex:auts-of-free-gl-one}, all GL-rack endomorphisms of $R_1$ have the form $(m,n)\mapsto(m+k,n+\ell)$ for some $(k,\ell)\in\Z^2$, and all of these endomorphisms are in fact automorphisms of $R_1$. In particular, $\eta_{R_1}$ is an automorphism of $R_1$, say \[\eta_{R_1}(m,n)=(m+k,n+\ell).\]
	
	By commutativity and the proof of Proposition \ref{prop:free-gl-one}, $\eta_{R_2}$ sends the image of the generator $(0,0)$ of $R_1$ under $\phi_0$ to \[\u^{\ell}_2\pi_{(Y,t)}^{k}\phi_0(0,0).\] Therefore, $\eta_{R_2}$ is completely determined by these powers of $\pi_{(Y,t)}$ and $\u_2$. Since $R_2$ was an arbitrary GL-rack, it follows that $Z(\glr)\subseteq\langle \Pi,\u\rangle$, as desired. It also follows from Example \ref{ex:auts-of-free-gl-one} and Proposition \ref{rmk:generators} that $Z(\glr)\cong\Z^2$, which proves claim \ref{Z1}.
	
	To prove claim \ref{Z2}, observe that $\Pi$ fixes $\mathbf{1}_\glq$. It follows from Proposition \ref{rmk:generators} that $Z(\glq)=\langle\u\rangle\cong\Z$, as desired. 
	Now, consider the full subcategory of $\glr$ whose objects are Legendrian racks. Using a similar argument as before, one can show using Corollary \ref{cor:leg-quandles} and Proposition \ref{rmk:generators} that the center of this category is \[\langle \Pi,\u\mid \Pi \u=\u\Pi,\,\Pi=\u^{-2}\rangle=\langle \u\rangle\cong\Z.\] 
	This proves claim \ref{Z2}. 
	Similarly, claim \ref{Z3} follows from Corollary \ref{cor:leg-quandles} and Proposition \ref{rmk:generators}.
\end{proof}
	
	\subsection{Categorical equivalence of racks and GL-quandles}
	Next, we show that the categories of racks and GL-quandles are isomorphic in a way that preserves mediality. 
	
	In a subsequent paper \cite{ta}, the author similarly showed that the category of racks is isomorphic to the category of Legendrian racks.
	
	\subsubsection{Construction of $F$}
	
	We begin by defining a functor $\F:\Rc\to\glq$. First, we define how $\F$ acts on objects.

\begin{prop}\label{prop:f}
	Given a rack $R=(X,s)$, define $\Fs:X\to S_X$ by \[x\mapsto \Fs_x:= \pi_R\inv s_x.\] Then $F(R):=(X,\Fs,\pi_R)$ is a GL-quandle.
\end{prop} 

\begin{proof}
	First, we show that $(X,\Fs)$ is a quandle. Part \ref{i} of Proposition \ref{prop:rack-center} and part \ref{iii} of Proposition \ref{lem:theta-inv} imply that, for all elements $x,y\in X$,
	\[
	\Fs_x\Fs_y=\pi\inv_Rs_x\pi\inv_R s_y=\pi^{-2}_R s_xs_y=\pi^{-2}_R s_{s_x(y)}s_x=\pi\inv_R s_{s_x(y)}\pi\inv_R s_x =\pi\inv_Rs_{\pi\inv_Rs_x(y)} \Fs_x = \Fs_{\Fs_x(y)}\Fs_x,
	\]
	so $(X,\Fs)$ is a rack.
	Moreover, part \ref{a2} of Proposition \ref{lem:theta-inv} and part \ref{i} of Proposition \ref{prop:rack-center} imply that
	\[
	\Fs_x(x)=\pi\inv_R s_x (x)= s_x\pi\inv_R (x)=s_xs_x\inv(x)=x,
	\]
	so $(X,\Fs)$ is a quandle.
	
	Next, we show that $\pi_R$ is a GL-structure on $(X,\Fs)$. Indeed, part \ref{i} of Proposition \ref{prop:rack-center} and part \ref{iii} of Proposition \ref{lem:theta-inv} imply that, for all $x\in X$,
	\begin{equation}\label{eq:cty-1}
		\pi_R \Fs_x = \pi_R \pi_R\inv s_x=\pi_R\inv s_x\pi_R=\pi\inv_R s_{\pi_R(x)}\pi_R=\Fs_{\pi_R(x)} \pi_R,
	\end{equation}
	so $\pi_R$ is a rack endomorphism of $(X,\Fs)$. Since $\pi_R$ is a bijection, we have $\pi_R\in\Aut (X,\Fs)$, as desired. Moreover, the third expression of equation (\ref{eq:cty-1}) equals $\Fs_x\pi_R$, so $\pi_R$ is a GL-structure.
\end{proof}

We now define how $\F$ acts on morphisms.

\begin{prop}
	For all racks $R=(X,s)$ and $S=(Y,t)$, and for all rack homomorphisms $f\in\Hom_\Rc(R,S)$, we have ${f\in \Hom_\glq(\F(R),\F(S))}$. So, if we define $F$ to fix $f$ as a set map, then $F$ is a covariant functor from $\Rc$ to $\glq$.
\end{prop}

\begin{proof}
	Certainly, $\F$ preserves the identity morphism and composition of morphisms, so we only need to verify that $f$ is a GL-rack homomorphism from $\F(R)=(X,\Fs,\pi_R)$ to $\F(S)=(Y,\Ft,\pi_S)$. Indeed, since $f\in \Hom_\Rc(R,S)$ and $\Pi\inv\in Z(\Rc)$, 
	\[
	f \Fs_x = f \pi_R\inv s_x = \pi_S\inv  f s_x=\pi_S\inv t_{f(x)} f=\Ft_{f(x)}  f
	\]
	for all $x\in X$, so $f\in\Hom_\Rc((X,\Fs),(Y,\Ft))$. Moreover, $f\pi_R=\pi_S f$ since $\Pi\in Z(\Rc)$, so $f$ is a GL-rack homomorphism.
\end{proof}

\subsubsection{Construction of $G$}
We now define a functor $\G:\glq\to\Rc$ as the restriction of a functor $\widetilde{\G}:\glr\to\Rc$ to $\glq$. First, we define how $\widetilde{G}$ acts on objects. 

\begin{prop}
	Given a GL-rack $R=(X,s,\u)$, define $\Gs:X\to S_X$ by \[x\mapsto\Gs_x:= \u s_x.\] Then $\widetilde{\G}(R):=(X,\Gs)$ is a rack.
\end{prop}

\begin{proof}
	Fix $x,y\in X$. Since $\u$ is a GL-structure and $(X,s)$ is a rack, 
	\[
	\Gs_{\Gs_x(y)} \Gs_x =\u s_{\u s_x(y)} \u s_x=\u^2 s_{s_x(y)} s_x=\u^2 s_xs_y=\u s_x\u s_y=\Gs_x\Gs_y,
	\]
	so $(X,\Gs)$ is a rack.
\end{proof}

Next, we define how $\widetilde{\G}$ acts on morphisms.

\begin{prop}
	For all GL-racks $R_1=(X,s,\u_1)$ and $R_2=(Y,t,\u_2)$ and GL-rack homomorphisms $g\in \Hom_\glr(R_1,R_2)$, we have $g\in \Hom_\Rc(\widetilde{\G}(R_1),\widetilde{\G}(R_2))$. So, if we define $\widetilde{G}$ to fix $g$ as a set map, then $\widetilde{\G}$ is a covariant functor from $\glr$ to $\Rc$, and $G$ is a functor from $\glq$ to $\Rc$.
\end{prop}

\begin{proof}
	Certainly, $\widetilde{\G}$ preserves the identity morphism and composition of morphisms, so we only need to verify that $g\in\Hom_\Rc((X,\Gs),(Y,\Gt))$. Indeed, since $g$ is a GL-rack homomorphism,
	\[
	g\Gs_x=g\u_1 s_x = \u_2 g s_x = \u_2 t_{g(x)} g = \Gt_{g(x)}g
	\]
	for all $x\in X$,
	as desired.
\end{proof}

\subsubsection{Isomorphism of categories}

Having defined $\F$ and $\G$, we are now ready to prove the main results of this section.

\begin{thm}\label{thm:isom}
	The functors $\F$ and $\G$ are isomorphisms of categories $\Rc\cong \glq$, and they restrict to isomorphisms $\Rc\med\cong\glq\med$.
\end{thm}

\begin{proof}
	To show that $\F$ and $\G$ are isomorphisms of categories, we only need to show that $\G\F$ and $\F\G$ fix the objects in the appropriate categories; $\G\F$ and $\F\G$ clearly fix the morphisms.
	To that end, let $R=(X,s)$ be a rack. To see that ${\G\F(R)=R}$, note that $\G\F(R)=(X,\hat{\Fs})$, where
	\[
	\hat{\Fs}_x=\pi_R \Fs_x = \pi_R \pi\inv_R s_x=s_x
	\]
	for all $x\in X$. That is, $\hat{\Fs}=s$, so $\G\F=\mathbf{1}_\Rc$, as desired. 
	
	Next, let $Q=(X,s,\u)$ be a GL-quandle. We must show that $FG(Q)=Q$. Note that, for all $x\in X$, we have \[\hat{s}_x\inv=s_x\inv\u\inv=\u\inv s_x\inv.\] Since $(X,s)$ is a quandle, we also have $x=s_x(x)$, so $s_x\inv(x)=x$. Now, to see that $FG(Q)=Q$, write $FG(Q)=(X,\widetilde{\Gs},\pi_{G(Q)})$. Part \ref{a2} of Proposition \ref{lem:theta-inv} implies that, for all elements $x,y\in X$,
	\[
	\widetilde{\Gs}_y(x)=\pi\inv_{G(Q)}\Gs_y(x)=\Gs_y\pi\inv_{G(Q)}(x)=\u s_y\hat{s}\inv_x(x)=s_y \u \u\inv s_x\inv(x)=s_y(x).
	\]
	Thus, $\widetilde{\Gs}_y=s_y$. Since $y\in X$ was arbitrary, this shows that $\widetilde{\Gs}=s$, as desired. Similarly,  \[\pi_{G(Q)}(x)=\hat{s}_x(x)=\u s_x(x)=\u (x)\] for all $x\in X$, so $\pi_{G(Q)}=\u$. 
	Hence, $FG=\mathbf{1}_\glq$, so $F$ and $G$ are isomorphisms of categories, as desired. 
	Since $\pi\inv$ and $\u$ are always rack automorphisms, the final claim follows straightforwardly from the definition of mediality that uses homomorphisms. 
\end{proof}
	
	\begin{cor}
		In the category of algebraic theories, the theory of racks and the theory of GL-quandles are isomorphic.
	\end{cor}
	
	\begin{proof}
		Since $F$ and $G$ are left and right adjoints to each other, they both preserve limits and colimits. In particular, they preserve finite products and filtered colimits, so by \cite{alg-theories}*{Lem.\ 3.8.3}, $F$ and $G$ are algebraic functors. Since $F$ and $G$ are equivalences of categories, the claim follows directly from \cite{alg-theories}*{Prop.\ 3.12.1}.
	\end{proof}
	
	\section{Tensor products of racks and GL-racks}\label{sec:tensors}
	In 2014, Crans and Nelson \cite{Hom}*{Sec.\ 8.1} categorified the results of Corollary \ref{cor:hom-racks} by considering universal-algebraic tensor products of medial quandles. However, these tensors remain unexplored in the literature. 
	
	With this motivation, we consider universal-algebraic tensor products of racks and GL-racks. We show that, unlike with groups, $\Rc$ and $\glr$ have tensor units. This suggests that $\glr$ and $\glr\med$ make natural settings for functorial invariants of Legendrian links.
	
	\subsection{Construction and universal property}
	
	We begin by considering tensor products of racks and GL-racks constructed via universal algebra. We discuss these tensors' universal properties and show that the induced symmetric monoidal structures on $\Rc\med$ and $\gla$ are closed.
	
	\begin{definition}\label{def:tensor}
		If $R_1=(X,r,\u_1)$ and $R_2=(Y,t,\u_2)$ are GL-racks, then we define their \emph{tensor product}, denoted by $R_1\otimes R_2$, to be the free GL-rack $\Free_\glr(X\times Y)$ modulo the smallest congruence such that the following hold for all $x,x_1,x_2\in X$ and $y,y_1,y_2\in Y$, writing $x\otimes y:=(x,y)$:
		\begin{enumerate}[(T1)]
			\item\label{T1} $s_{x\otimes y_2}(x\otimes y_1)\sim x\otimes t_{y_2}(y_1)$,\quad and \quad$s\inv_{x\otimes y_2}(x\otimes y_1)\sim x\otimes t\inv_{y_2}(y_1)$.
			\item\label{T2} $s_{x_2 \otimes y}(x_1\otimes y)\sim r_{x_2}(x_1)\otimes y$,\quad and \quad$s\inv_{x_2 \otimes y}(x_1\otimes y)\sim r\inv_{x_2}(x_1)\otimes y$.
			\item\label{T3} $\u(x\otimes y)\sim \u_1(x)\otimes y\sim x\otimes \u_2(y).$ 
		\end{enumerate}
		We also define the \emph{medial tensor product} of $R_1$ and $R_2$ to be \[R_1\otimes\med R_2:=(R_1\otimes R_2)\med.\] 
		Similarly, we define the tensor product of two racks $R_1=(X,r)$ and $R_2=(Y,t)$ as $\Free_\Rc(X\times Y)$ modulo the smallest congruence generated by relations \ref{T1} and \ref{T2} above, and we define the medial tensor product of racks similarly.
	\end{definition}
	
	\begin{remark}\label{tensor-u-inv}
		Relation \ref{T3} implies a similar relation involving $\u\inv$. Namely, in $(X,r,\u_1)\otimes (Y,t,\u_2)$,  \[\u\inv(x\otimes y)=\u_1\inv (x)\otimes y=x\otimes\u\inv_2(y)\] for all $x\in X$ and $y\in Y$. To see this, compute \[\u_1\inv(x)\otimes y=\u\inv\u(\u_1\inv(x)\otimes y)=\u\inv(\u_1\u_1\inv(x)\otimes y)=\u\inv(x\otimes y),\] and similarly for $x\otimes\u\inv_2(y)$.
	\end{remark}
	
	We note that the medialized associated quandles of the tensor products of racks in Definition \ref{def:tensor} recover the tensor products of medial quandles that Crans and Nelson \cite{Hom}*{Sec.\ 8.1} introduced in 2014. On the other hand, these tensors are distinct from the tensor products that Kamada \cite{tensors}*{Def.\ 3.1} introduced in 2021, which are sets without canonical rack structures. 
	
	\subsubsection{Bihomomorphisms} The discussion in \cite{alg-theories}*{p.\ 171} shows that the tensor products in Definition \ref{def:tensor} are each characterized by a universal factorizing property. Before we can state this property, we will need the following definition from universal algebra. This is a special case of a more general definition for arbitrary algebraic theories; see, for example, \cite{alg-theories}*{Def.\ 3.10.2}.
	
	\begin{definition}
	Let $(X,r)$, $(Y,s),$ and $(Z,t)$ be racks. We say that a map $\beta_0:X\times Y\to Z$ is a \emph{rack bihomomorphism} if, for all $x\in X$ and $y\in Y$, the restricted maps $\beta_0(-,y):X\to Z$ and $\beta_0(x,-):Y\to Z$ are rack homomorphisms. \emph{GL-rack bihomomorphisms} are defined similarly.
\end{definition}

\begin{example}\label{ex:bihoms}
	Let $R=(X,r)$ be a rack, and let $S=(Y,s)$ and $T=(Z,t)$ be medial racks. Recall from Corollary \ref{cor:hom-racks} that $H_1:=\Hom_\Rc(R,S)$, $H_2:=\Hom_\Rc(S,T)$, and $H_3:=\Hom_\Rc(R,T)$ each have a canonical rack structure.
	
	In analogy with the composition of $A$-linear maps in the category of modules over a ring $A$, the composition map $\beta_0:H_1\times H_2\to H_3$ is a rack bihomomorphism. To see this, fix a homomorphism $g\in H_2$. For all homomorphisms $\phi,\psi\in H_1$ and elements $x\in X$, we have
	\[
	\beta_0(\Fs_\psi(\phi),g)(x)=g\Fs_\psi(\phi)(x)=gs_{\psi(x)}\phi(x)=t_{g\psi(x)}g\phi(x)=\widetilde{t}_{g\psi} g\phi(x)=\widetilde{t}_{\beta_0(\psi,g)}\beta_0(\phi,g)(x),
	\]
	so the restriction $\beta_0(-,g):H_1\to H_3$ is a rack homomorphism. Similarly, for all homomorphisms $f\in H_1$, the restriction $\beta_0(f,-):H_2\to H_3$ is a homomorphism, so $\beta_0$ is a bihomomorphism. 
	
	It is straightforward to verify using Theorem \ref{thm:homsets} that if $R$, $S$, and $T$ also have GL-structures, then the restriction of $\beta_0$ to the respective hom-sets in $\glr$ is a GL-rack bihomomorphism.
\end{example}

	\subsubsection{Universal property of tensor products}
	As suggested in Example \ref{ex:bihoms}, the definition of a rack bihomomorphism analogizes the definition of a bilinear map in the category of modules over a ring. The universal property of tensor products, which we state below, extends this analogy.

\begin{prop}\label{prop:univ-tensor}
	Let $R_1$ and $R_2$ be racks. Then $R_1\otimes R_2$ is characterized up to isomorphism by a universal property. Namely, there exists a rack bihomomorphism $\psi$ from $R_1\times R_2$ to $R_1\otimes R_2$ such that, for all racks $R_3$ and rack bihomomorphisms $\beta_0:R_1\times R_2\to R_3$, there exists a unique rack homomorphism $\beta:R_1\otimes R_2\to R_3$ such that $\beta_0=\beta \psi$. In particular, the following diagram commutes:
	\begin{equation}\label{diag:tensors}\begin{tikzcd}
			{R_1\times R_2} & {R_1\otimes R_2} \\
			& R_3
			\arrow["\psi", from=1-1, to=1-2]
			\arrow["{\beta_0}"', from=1-1, to=2-2]
			\arrow["\beta", dashed, from=1-2, to=2-2]
	\end{tikzcd}\end{equation}
	A similar result holds for GL-racks.
\end{prop}

\begin{proof}
	This is a consequence of universal algebra. Regard racks and GL-racks as algebraic theories with binary operations $s_-(-)$ and $s\inv_-(-)$, along with unary operations $\u$ and $\u\inv$ for GL-racks. Due to Remark \ref{tensor-u-inv}, the tensor product in Definition \ref{def:tensor} is precisely the tensor product constructed in the proof of \cite{alg-theories}*{Thm.\ 3.10.3}, and $\psi(x,y)=x\otimes y$ for all $(x,y)\in X\times Y$. Thus, the claim follows from the discussion in \cite{alg-theories}*{p.\ 171}.
\end{proof}
	
	\begin{example}\label{ex:left-dist}
	A rack $(X,s)$ is called \emph{left-distributive} if \begin{equation}\label{eq:left-dist}
		s_{s_a(b)}(x)=s_{s_a(x)} (s_b(x))
	\end{equation}
	for all $a,b,x\in X$. The name comes from the fact that, in terms of the right-distributive binary operation $\tr$ often used in the literature, equation (\ref{eq:left-dist}) states that \[x\tr(b\tr a)=(x\tr b)\tr (x\tr a).\] 
	
	For example, it is straightforward to verify that medial quandles are left-distributive. 
	Moreover, self-distributive quasigroups are precisely left-distributive quandles satisfying an axiom called the \emph{Latin} condition; see, for example, \cite{quandlebook}*{p.\ 143}.
	
	Evidently, a rack $R=(X,s)$ is left-distributive if and only if, for all $x\in X$, the map $X\to X$ defined by $y\mapsto s_y(x)$ is an endomorphism of $R$. Equivalently, the map $\beta_0:X\times X\to X$ defined by $(x,y)\mapsto s_y(x)$ is a bihomomorphism from $R\times R$ to $R$. 
	In this case, the universal property of $R\otimes R$ implies the existence of a unique homomorphism $\beta\in\Hom_\Rc(R\otimes R, R)$ such that $\beta(x\otimes y)=s_y(x)$ for all $x,y\in X$. 
\end{example}

\subsubsection{Internal hom-tensor adjunctions}
Mirroring the work of Crans and Nelson \cite{Hom}*{Thm.\ 12} on the category of medial quandles, the next result continues the analogy with the category of modules over a ring by describing internal hom-tensor adjunctions in $\Rc\med$ and $\glr\med$.

\begin{thm}\label{thm:scmc}
	The categories $\Rc\med$ and $\gla$ are closed symmetric monoidal with respect to the medial tensor product $\otimes\med$ in each category and the closed structures $\Hom_{\Rc\med}(-,-)$ and $\Hom_{\gla}(-,-)$ from Corollary \ref{cor:hom-racks} and Theorem \ref{thm:homsets}, respectively.
\end{thm}

\begin{proof}
	Recall that the algebraic theories of medial racks and medial GL-racks are commutative. Therefore, both claims are special cases of a general result for commutative algebraic theories; see, for example, \cite{alg-theories}*{Thm.\ 3.10.3}. 
	
	In particular, $\otimes\med$ is precisely the tensor product constructed in the proof of \cite{alg-theories}*{Thm.\ 3.10.3}, and the tensor unit is the free rack (resp.\ GL-rack) $L$ on one element. Since permutation racks are medial, Example \ref{ex:free-one} (resp.\ Proposition \ref{prop:free-gl-one}) shows that $L$ is medial.
\end{proof}

\subsection{Nonmedial tensor products}
One may ask how much of the structure in Theorem \ref{thm:scmc} remains if we drop the mediality assumption. By \cite{alg-theories}*{Thm.\ 3.10.3}, we lose the closed structure because racks and GL-racks are noncommutative algebraic theories. Nevertheless, we show that tensor products in $\Rc$ and $\glr$ surprisingly retain tensor units.

Although universal-algebraic tensor products with the appropriate universal properties exist for all algebraic theories, they are often not well-behaved in noncommutative algebraic theories. For example, the universal-algebraic tensor product of groups $G\otimes H$ is isomorphic to the usual $\Z$-module tensor product $G^{\operatorname{ab}}\otimes_\Z H^{\operatorname{ab}}$ of the abelianizations of $G$ and $H$; see, for example, \cite{alg-theories}*{p.\ 171}. In particular, no tensor unit exists. 

	In this subsection, we show that the pathologies in the previous paragraph do \emph{not} apply to tensor products of racks or GL-racks, even though racks and GL-racks are noncommutative algebraic theories. In particular, while the universal-algebraic tensor product of groups is always abelian, tensor products of racks and GL-racks are not necessarily medial.



\begin{lemma}\label{lem:theta-in-tensor}
	Let $R_1=(X,r)$, $R_2=(Y,s)$, and $R_3=(Z,t)$ be racks, and let $\beta_0:X\times Y\to Z$ be a rack bihomomorphism from $R_1\times R_2$ to $R_3$. Then 
	for all integers $k\in\Z$, elements $x\in X$, and elements $y\in Y$,  \[\beta_0(\pi_{R_1}^k(x),y)=\beta_0(x,\pi_{R_2}^k(y)).\]
\end{lemma}

\begin{proof}
		For all $y\in Y$, the restriction $\beta_0(-,y):X\to Z$ is a homomorphism from $R_1$ to $R_3$. For all integers $k\in\Z$, the inclusion $\Pi^k\in Z(\Rc)$ implies that  \[\beta_0(\pi_{R_1}^k(x),y)=\pi^k_{R_3}\beta_0(x,y)\] for all $x\in X$. A similar argument using the restriction $\beta_0:(x,-):Y\to Z$ yields the claim. 
\end{proof}

\begin{thm}\label{prop:tensor-unit}
	Let $F$ be the free rack on one element. Then for all racks $R$, we have natural isomorphisms $R\otimes F\cong R\cong F\otimes R$.
	
	Similarly, let $L$ be the free GL-rack on one element. Then for all GL-racks $R$, we have natural isomorphisms $R\otimes L\cong R\cong L\otimes R$.
\end{thm}

	\begin{proof}
	First, let $R=(X,s)$ be a rack. We will show that $R\otimes F\cong R$; the proof that $F\otimes R\cong R$ is similar, and naturality is straightfoward to verify from there. Identify $F=(\Z,\sigma)\perm$ 
	as in Example \ref{ex:free-one}.
	
	By Proposition \ref{prop:univ-tensor}, it will suffice to show that $R$ satisfies the universal property of $R\otimes F$ with $\psi:X\times \Z\to X$ defined by \[(x,k)\mapsto \pi^k_R(x).\] First, we show that $\psi$ is a rack bihomomorphism. Fix $x\in X$, and let $s^\Z:\Z\to \{\sigma\}$ denote the rack structure of $F$. For all $k,n\in \Z$,
	\[
	\psi(x,s^\Z_k(n))=\psi(x,\sigma(n))=\psi(x,n+1)=\pi^{n+1}_R(x)=\pi^n_R s_x(x)=s_{\pi_R^k(x)}\pi_R^n(x)=s_{\psi(x,k)}\psi(x,n).
	\]
	In the fifth equality, we have used part \ref{i} of Proposition \ref{prop:rack-center} and part \ref{iii} of Proposition \ref{lem:theta-inv}. So, ${\psi(x,-):\Z\to Y}$ is a rack homomorphism from $F$ to $R'$.
	Next, fix $k\in\Z$. For all $x,y\in X$,
	\[
	\psi(s_x(y),k)=\pi^k_Rs_x(y)=s_x\pi^k_R(y)=s_{\pi_R^k(x)}\pi^k_R(y)=s_{\psi(x,k)}\psi(y,k)
	\]
	by part \ref{i} of Proposition \ref{prop:rack-center} and part \ref{iii} of Proposition \ref{lem:theta-inv}, so $\psi(-,k):X\to Y$ is a rack homomorphism from $R$ to $R'$. Therefore, $\psi$ is a rack bihomomorphism, as desired.
	
	Now, let $R'=(Y,t)$ be a rack, and let $\beta_0:X\times \Z\to Y$ be a rack bihomomorphism from $R\times F$ to $R'$. Define $\beta:X\to Y$ by $x\mapsto \beta_0(x,0)$. Then $\beta$ is precisely the restriction $\beta_0(-,0):X\to Y$, so $\beta$ is a rack homomorphism, as desired. Since $\pi_F=\sigma$, Lemma \ref{lem:theta-in-tensor} implies that
	\[
	\beta\psi(x,k)=\beta\pi^k_R(x)=\beta_0(\pi^k_R(x),0)=\beta_0(x,\sigma^k(0))=\beta_0(x,k)
	\]
	for all $(x,k)\in X\times\Z$, so diagram (\ref{diag:tensors}) commutes. Finally, uniqueness follows from the surjectivity of $\psi$. Hence, $R$ satisfies the universal property of $R\otimes F$, so $R\cong R\otimes F$, as claimed.
	
	For the second part of the claim, let $R=(X,s,\u)$ be a GL-rack. Identify $L=(\Z^2,\sigma,\u_0)\perm$ as in Proposition \ref{prop:free-gl-one}. Again, it suffices to show that $R$ satisfies the universal property of $R\otimes L$ with $\psi:X\times \Z^2$ defined by \[(x,m,n)\mapsto \u^n \pi^m_R(x).\] The proof that $\psi$ is a GL-rack bihomomorphism is similar to the argument given above; we leave the details to the reader.
	
	Given a GL-rack $R'=(Y,t,\u_2)$ and a GL-rack bihomomorphism $\beta_0:X\times\Z^2\to Y$ from $R\times L$ to $R'$, define ${\beta:X\to Y}$ by $x\mapsto \beta_0(x,0,0)$. Once again, $\beta$ is a GL-rack homomorphism since it is the restriction $\beta_0(-,0,0):X\to Y$. The commutativity of diagram (\ref{diag:tensors}) and uniqueness of $\beta$ are shown in a similar way as before, so $R\cong R\otimes L$.
\end{proof}
	
	Recall that the free rack on one element is medial, and similarly for the free GL-rack on one element. Therefore, Theorem \ref{prop:tensor-unit} implies the following result---a surprising statement considering that the universal-algebraic tensor product of groups is necessarily abelian.

	\begin{cor}
		Even if one of the tensor factors is medial, tensor products of racks are not necessarily medial, and similarly for tensor products of GL-racks.
	\end{cor}
	
	It would be interesting to study more properties and applications of the tensor products in Definition \ref{def:tensor}.
	In Section \ref{sec:questions}, we propose future work in this direction.

	\section{Directions for future work}\label{sec:questions}
	We conclude by proposing questions for further research. In the following, let $\GG(\Lambda)$ denote the fundamental GL-rack of a Legendrian link $\Lambda$ (see \cite{karmakar}).
	
	\begin{enumerate}
		\item Use GL-racks to distinguish the conjecturally nonequivalent Legendrian knots listed in the Legendrian knot atlases \citelist{\cite{atlas}\cite{atlas-site}\cite{new-atlas}}. We note that this may be a suitable problem for automation; as discussed in \citelist{\cite{atlas}\cite{us}\cite{new-atlas}}, grid diagrams and Legendrian knot mosaics could yield suitably discrete front projections for computer programs to traverse. 
		\item Classify more families of GL-racks using Theorem \ref{cor:number}.
		\item Compute more GL-rack automorphism groups using Proposition \ref{thm:aut-grp}.
		\item In light of Proposition \ref{prop:medial-enhance}, do there exist Legendrian knots $\Lambda_1$ and $\Lambda_2$ that have the same coloring numbers (see \citelist{\cite{karmakar}\cite{bi}\cite{ceniceros}}) by a medial GL-rack $M$ such that $\Hom_\glr(\GLR(\Lambda_1),M)$ and $\Hom_\glr(\GLR(\Lambda_2),M)$ are nonisomorphic as medial GL-racks? Also, do there exist such $\Lambda_1$ and $\Lambda_2$ whose classical invariants are identical? A positive answer would show that $\Hom_\glr(\GG(\Lambda),M)$ is a proper enhancement of the coloring number as an invariant of $\Lambda$.
		\item Let $\mathbb{F}$ be a field, and let $M$ be a medial GL-rack. In 2023, Elhamdadi, Nunez, and Singh \cite{enhancements}*{Thms.\ 4.2 and 5.1} properly enhanced medial quandle-valued invariants of smooth links using $\mathbb{F}$-algebra homomorphisms between quandle rings and colorings of smooth links by idempotents of quandle rings. Do similar proper enhancements of $\Hom_\glr(\GLR(\Lambda),M)$ exist? 
		\item In light of Example \ref{ex:left-dist}, do tensor products of left-distributive racks, medial quandles (cf.\ \cite{Hom}*{Sec.\ 8.1}), or Latin quandles have any interesting properties?
		\item Theorem \ref{prop:tensor-unit} shows that universal-algebraic tensor products make $\Rc$ and $\glr$ into \emph{symmetric magmoidal categories with units}, leading us to ask what their \emph{unital nuclei} are; see \cite{magmoidal}*{Sec.\ 2.5}. In particular, is $\otimes$ associative? If so, then $\otimes$ induces symmetric monoidal structures on $\Rc$ and $\glr$.
			\item Can tensor products of racks or GL-racks be used to define new invariants of smooth links or Legendrian links?
		\item Theorem \ref{thm:scmc} implies that $\gla$ enriches over itself.
		In this light, what are the applications of enriched category theory to medial GL-racks and invariants of Legendrian links? 
		\item Study 4-Legendrian biracks, which Kimura \cite{leg-cocycle}*{Sec.\ 4} introduced in 2024.
		\item In light of Remark \ref{rmk:virtual}, can our results about GL-racks be generalized to virtual racks or virtual biracks?
\item In 2015, Cahn and Levi \cite{virtual-leg} introduced \emph{virtual Legendrian knots}, which are Legendrian knots in the spherical cotangent bundle of a surface equipped with the natural contact structure. Can GL-racks be used to define invariants of virtual Legendrian links (cf.\ Remark \ref{rmk:virtual})?
		\item \emph{Transverse knots} are knots that lie everywhere transverse to the standard contact structure on $\R^3$; see \cite{etnyre}*{Sec.\ 2.4}. Can one define rack-theoretic invariants of transverse knots?
		\item The fundamental quandle of a smooth link is interpreted topologically as the set of homotopy classes of paths from a basepoint in the link complement to the boundary of the link complement with several restrictions; see, for example, \cite{quandlebook}*{p.\ 125}. Is there a similar contact-topological interpretation of the GL-rack of a Legendrian link?
\item Are there formulas for the number of GL-racks and medial GL-racks of a given finite order?
	\end{enumerate}
	
	\section*{Declarations}
	
	\begin{competing}
		The author declares that he has no conflict of interest.
	\end{competing}
		
			\bibliography{references}  

@article {ceniceros,
	AUTHOR = {Ceniceros, Jose and Elhamdadi, Mohamed and Nelson, Sam},
	TITLE = {Legendrian rack invariants of {L}egendrian knots},
	JOURNAL = {Commun. Korean Math. Soc.},
	FJOURNAL = {Korean Mathematical Society. Communications},
	VOLUME = {36},
	YEAR = {2021},
	NUMBER = {3},
	PAGES = {623--639},
	ISSN = {1225-1763,2234-3024},
	MRCLASS = {57K12 (57K10)},
	MRNUMBER = {4292403},
	MRREVIEWER = {Pedro\ Lopes},
	DOI = {10.4134/CKMS.c200251},
	URL = {https://doi.org/10.4134/CKMS.c200251},
}

@misc{original,
	title={On Rack Invariants Of {L}egendrian Knots}, 
	author={Dheeraj Kulkarni and T. V. H. Prathamesh},
	year={2017},
	eprint={1706.07626},
	archivePrefix={arXiv},
	primaryClass={math.GT},
	url={https://arxiv.org/abs/1706.07626}, 
	NOTE = {Preprint, arXiv:1706.07626 [math.GT].},
}

@article {karmakar,
	AUTHOR = {Karmakar, Biswadeep and Saraf, Deepanshi and Singh, Mahender},
	TITLE = {Generalized {L}egendrian racks of {L}egendrian links},
	JOURNAL = {J. Topol. Anal.},
	FJOURNAL = {Journal of Topology and Analysis},
	VOLUME = {18},
	YEAR = {2026},
	NUMBER = {6},
	PAGES = {1857--1874},
	ISSN = {1793-5253,1793-7167},
	MRCLASS = {57K33 (57K10 57K12)},
	MRNUMBER = {5051798},
	DOI = {10.1142/S1793525326500020},

}

@article {karmakar2,
	AUTHOR = {Karmakar, Biswadeep and Saraf, Deepanshi and Singh, Mahender},
	TITLE = {Beck modules over generalized {L}egendrian racks},
	JOURNAL = {Internat. J. Algebra Comput.},
	FJOURNAL = {International Journal of Algebra and Computation},
	VOLUME = {36},
	YEAR = {2026},
	NUMBER = {4},
	PAGES = {407--424},
	ISSN = {0218-1967,1793-6500},
	MRCLASS = {57K12 (18N40 20N02 57K33)},
	MRNUMBER = {5077711},
	DOI = {10.1142/S0218196726500153},
}

@article {joyce,
	AUTHOR = {Joyce, David},
	TITLE = {A classifying invariant of knots, the knot quandle},
	JOURNAL = {J. Pure Appl. Algebra},
	FJOURNAL = {Journal of Pure and Applied Algebra},
	VOLUME = {23},
	YEAR = {1982},
	NUMBER = {1},
	PAGES = {37--65},
	ISSN = {0022-4049,1873-1376},
	MRCLASS = {57M25 (20F29 20N05 53C35)},
	MRNUMBER = {638121},
	MRREVIEWER = {Mark\ E.\ Kidwell},
	DOI = {10.1016/0022-4049(82)90077-9},
	URL = {https://doi.org/10.1016/0022-4049(82)90077-9},
}

@article {bi,
	AUTHOR = {Kimura, Naoki},
	TITLE = {Bi-{L}egendrian rack colorings of {L}egendrian knots},
	JOURNAL = {J. Knot Theory Ramifications},
	FJOURNAL = {Journal of Knot Theory and its Ramifications},
	VOLUME = {32},
	YEAR = {2023},
	NUMBER = {4},
	PAGES = {Paper No. 2350029, 16},
	ISSN = {0218-2165,1793-6527},
	MRCLASS = {57K12 (57K10 57K33)},
	MRNUMBER = {4586264},
	MRREVIEWER = {Mohamed\ Elhamdadi},
	DOI = {10.1142/S0218216523500293},
	URL = {https://doi.org/10.1142/S0218216523500293},
}

@article {us,
	AUTHOR = {Kipe, Margaret and Pezzimenti, Samantha and Schaumann, Leif
	and Ta, L{\d\uhorn}c and Wong, Tony W. H.},
	TITLE = {Bounds on the mosaic number of {L}egendrian knots},
	JOURNAL = {J. Knot Theory Ramifications},
	FJOURNAL = {Journal of Knot Theory and its Ramifications},
	VOLUME = {34},
	YEAR = {2025},
	NUMBER = {12},
	PAGES = {Paper No. 2550055, 52},
	ISSN = {0218-2165,1793-6527},
	MRCLASS = {57K10 (57K33)},
	MRNUMBER = {4966653},
	DOI = {10.1142/S0218216525500555},
}

@incollection {etnyre,
	AUTHOR = {Etnyre, John B.},
	TITLE = {Legendrian and transversal knots},
	BOOKTITLE = {Handbook of knot theory},
	PAGES = {105--185},
	PUBLISHER = {Elsevier B. V., Amsterdam},
	YEAR = {2005},
	ISBN = {0-444-51452-X},
	MRCLASS = {57R17 (53D35 57M25 57M27)},
	MRNUMBER = {2179261},
	MRREVIEWER = {Lenhard\ L.\ Ng},
	DOI = {10.1016/B978-044451452-3/50004-6},
	URL = {https://doi.org/10.1016/B978-044451452-3/50004-6},
}

@article {aut,
	AUTHOR = {Elhamdadi, Mohamed and Macquarrie, Jennifer and Restrepo,
	Ricardo},
	TITLE = {Automorphism groups of quandles},
	JOURNAL = {J. Algebra Appl.},
	FJOURNAL = {Journal of Algebra and its Applications},
	VOLUME = {11},
	YEAR = {2012},
	NUMBER = {1},
	PAGES = {1250008, 9},
	ISSN = {0219-4988,1793-6829},
	MRCLASS = {20B25 (20N02)},
	MRNUMBER = {2900878},
	MRREVIEWER = {Alexander\ Zvonkin},
	DOI = {10.1142/S0219498812500089},
	URL = {https://doi.org/10.1142/S0219498812500089},
}

@article {intro,
	AUTHOR = {Nelson, Sam},
	TITLE = {What is {$\ldots$} a quandle?},
	JOURNAL = {Notices Amer. Math. Soc.},
	FJOURNAL = {Notices of the American Mathematical Society},
	VOLUME = {63},
	YEAR = {2016},
	NUMBER = {4},
	PAGES = {378--380},
	ISSN = {0002-9920,1088-9477},
	MRCLASS = {57M27 (20N99)},
	MRNUMBER = {3444659},
	DOI = {10.1090/noti1360},
	URL = {https://doi.org/10.1090/noti1360},
}

@article {survey,
	AUTHOR = {Elhamdadi, Mohamed},
	TITLE = {A survey of racks and quandles: {S}ome recent developments},
	JOURNAL = {Algebra Colloq.},
	FJOURNAL = {Algebra Colloquium},
	VOLUME = {27},
	YEAR = {2020},
	NUMBER = {3},
	PAGES = {509--522},
	ISSN = {1005-3867,0219-1733},
	MRCLASS = {20N02 (20C99 57K12)},
	MRNUMBER = {4141628},
	MRREVIEWER = {David\ Stanovsk\'y},
	DOI = {10.1142/S1005386720000425},
	URL = {https://doi.org/10.1142/S1005386720000425},
}

@book {book,
	AUTHOR = {Nosaka, Takefumi},
	TITLE = {Quandles and topological pairs},
	SERIES = {SpringerBriefs in Mathematics},
	SUBTITLE = {Symmetry, knots, and cohomology},
	PUBLISHER = {Springer, Singapore},
	YEAR = {2017},
	PAGES = {ix+136},
	ISBN = {978-981-10-6792-1; 978-981-10-6793-8},
	MRCLASS = {57M27 (20J06)},
	MRNUMBER = {3729413},
	MRREVIEWER = {Markus\ Szymik},
	DOI = {10.1007/978-981-10-6793-8},
	URL = {https://doi.org/10.1007/978-981-10-6793-8},
}

@article {Hom,
	AUTHOR = {Crans, Alissa S. and Nelson, Sam},
	TITLE = {Hom quandles},
	JOURNAL = {J. Knot Theory Ramifications},
	FJOURNAL = {Journal of Knot Theory and its Ramifications},
	VOLUME = {23},
	YEAR = {2014},
	NUMBER = {2},
	PAGES = {1450010, 18},
	ISSN = {0218-2165,1793-6527},
	MRCLASS = {57M27 (57M25)},
	MRNUMBER = {3197054},
	MRREVIEWER = {Bruno\ P.\ Zimmermann},
	DOI = {10.1142/S0218216514500102},
	URL = {https://doi.org/10.1142/S0218216514500102},
}

@book {leg-cocycle,
	AUTHOR = {Kimura, Naoki},
	TITLE = {Rack coloring invariants of {L}egendrian knots},
	NOTE = {Thesis (Ph.D.)--Waseda University Graduate School of Fundamental Science and Engineering},
	PUBLISHER = {},
	YEAR = {2024},
	PAGES = {},
	MRCLASS = {},
	MRNUMBER = {},
	URL = {https://waseda.repo.nii.ac.jp/record/2002429/files/Honbun-9489.pdf},
}

@article {atlas,
	AUTHOR = {Chongchitmate, Wutichai and Ng, Lenhard},
	TITLE = {An atlas of {L}egendrian knots},
	JOURNAL = {Exp. Math.},
	FJOURNAL = {Experimental Mathematics},
	VOLUME = {22},
	YEAR = {2013},
	NUMBER = {1},
	PAGES = {26--37},
	ISSN = {1058-6458,1944-950X},
	MRCLASS = {57M25 (53Cxx 57R17)},
	MRNUMBER = {3038780},
	DOI = {10.1080/10586458.2013.750221},
	URL = {https://doi.org/10.1080/10586458.2013.750221},
}

@misc{atlas-site,
	author = {Bhattacharyya, Nilangshu and Cox, Cyrus and Murray, Justin and Pandikkadan, Adithyan and Vela-Vick, Shea and Wu, Angela},
	title = {Legendrian Knot Atlas},
	howpublished = {\url{https://www.math.lsu.edu/~knotatlas/legendrian/index.html}},
	url      	= {https://www.math.lsu.edu/~knotatlas/legendrian/index.html},
	note        = {\url{https://www.math.lsu.edu/~knotatlas/legendrian/index.html}. Accessed: 2025-3-15},
	year = {n.d.},
}

@manual{GAP4,
	organization = "The GAP~Group",
	title        = "{GAP -- Groups, Algorithms, and Programming,
	Version 4.14.0}",
	year         = 2024,
	url          = "\url{https://www.gap-system.org}",
}

@article {library,
	AUTHOR = {Vojtěchovský, Petr and Yang, Seung Yeop},
	TITLE = {Enumeration of racks and quandles up to isomorphism},
	JOURNAL = {Math. Comp.},
	FJOURNAL = {Mathematics of Computation},
	VOLUME = {88},
	YEAR = {2019},
	NUMBER = {319},
	PAGES = {2523--2540},
	ISSN = {0025-5718,1088-6842},
	MRCLASS = {20N02 (16T25 57M27)},
	MRNUMBER = {3957904},
	MRREVIEWER = {David\ Stanovsk\'y},
	DOI = {10.1090/mcom/3409},
	URL = {https://doi.org/10.1090/mcom/3409},
}

@article {enhancements,
	AUTHOR = {Elhamdadi, Mohamed and Nunez, Brandon and Singh, Mahender},
	TITLE = {Enhancements of link colorings via idempotents of quandle
	rings},
	JOURNAL = {J. Pure Appl. Algebra},
	FJOURNAL = {Journal of Pure and Applied Algebra},
	VOLUME = {227},
	YEAR = {2023},
	NUMBER = {10},
	PAGES = {Paper No. 107400, 16},
	ISSN = {0022-4049,1873-1376},
	MRCLASS = {17D99 (16S34 20N02)},
	MRNUMBER = {4579329},
	MRREVIEWER = {Emanuele\ Zappala},
	DOI = {10.1016/j.jpaa.2023.107400},
	URL = {https://doi.org/10.1016/j.jpaa.2023.107400},
}

@misc{library-site, 
	AUTHOR = {Vojtěchovský, Petr and Yang, Seung Yeop}, 
	TITLE = "Racks and quandles of small orders",
	YEAR = "2018",
	HOWPUBLISHED = "\url{https://www.cs.du.edu/~petr/libraries_of_algebraic_structures.html}", 
	NOTE = {{\url{https://www.cs.du.edu/~petr/libraries_of_algebraic_structures.html}. Accessed: 2025-01-03}},
	URL = "https://www.cs.du.edu/~petr/libraries_of_algebraic_structures.html",
}

@book{quandlebook,
	AUTHOR = {Elhamdadi, Mohamed and Nelson, Sam},
	TITLE = {Quandles},
	SUBTITLE = {An introduction to the algebra of knots},
	SERIES = {Student Mathematical Library},
	VOLUME = {74},
	PUBLISHER = {American Mathematical Society, Providence, RI},
	YEAR = {2015},
	PAGES = {x+245},
	ISBN = {978-1-4704-2213-4},
	MRCLASS = {57M27 (57M25 57Q45)},
	MRNUMBER = {3379534},
	MRREVIEWER = {Frederick\ Norwood},
	DOI = {10.1090/stml/074},
	URL = {https://doi.org/10.1090/stml/074},
}

@article {takasaki,
	AUTHOR = {Takasaki, Mituhisa},
	TITLE = {Abstraction of symmetric transformations},
	JOURNAL = {T\^ohoku Math. J.},
	FJOURNAL = {The T\^ohoku Mathematical Journal},
	VOLUME = {49},
	YEAR = {1943},
	PAGES = {145--207},
	ISSN = {0040-8735,1881-2015},
	MRCLASS = {20.0X},
	MRNUMBER = {21002},
	MRREVIEWER = {S.\ Kakutani},
}

@article {fenn,
	AUTHOR = {Fenn, Roger and Rourke, Colin},
	TITLE = {Racks and links in codimension two},
	JOURNAL = {J. Knot Theory Ramifications},
	FJOURNAL = {Journal of Knot Theory and its Ramifications},
	VOLUME = {1},
	YEAR = {1992},
	NUMBER = {4},
	PAGES = {343--406},
	ISSN = {0218-2165,1793-6527},
	MRCLASS = {57M25 (57N10)},
	MRNUMBER = {1194995},
	DOI = {10.1142/S0218216592000203},
	URL = {https://doi.org/10.1142/S0218216592000203},
}

@article {rack-roll,
	AUTHOR = {Grøsfjeld, Tobias},
	TITLE = {Thesaurus racks: {C}ategorizing rack objects},
	JOURNAL = {J. Knot Theory Ramifications},
	FJOURNAL = {Journal of Knot Theory and its Ramifications},
	VOLUME = {30},
	YEAR = {2021},
	NUMBER = {4},
	PAGES = {Paper No. 2150019, 18},
	ISSN = {0218-2165,1793-6527},
	MRCLASS = {18C40 (16B50 20J15 57K12)},
	MRNUMBER = {4272643},
	MRREVIEWER = {Markus\ Szymik},
	DOI = {10.1142/S021821652150019X},
	URL = {https://doi.org/10.1142/S021821652150019X},
}

@book {alg-theories,
	AUTHOR = {Borceux, Francis},
	TITLE = {Handbook of categorical algebra, volume 2},
	SUBTITLE = {Categories and structures},
	SERIES = {Encyclopedia of Mathematics and its Applications},
	VOLUME = {51},
	NOTE = {},
	PUBLISHER = {Cambridge University Press, Cambridge},
	YEAR = {1994},
	PAGES = {xviii+443},
	ISBN = {0-521-44179-X},
	MRCLASS = {18-02 (18Exx)},
	MRNUMBER = {1313497},
	MRREVIEWER = {Martin\ Hyland},
}

@article {center,
	AUTHOR = {Szymik, Markus},
	TITLE = {Permutations, power operations, and the center of the category
	of racks},
	JOURNAL = {Comm. Algebra},
	FJOURNAL = {Communications in Algebra},
	VOLUME = {46},
	YEAR = {2018},
	NUMBER = {1},
	PAGES = {230--240},
	ISSN = {0092-7872,1532-4125},
	MRCLASS = {18C10 (20N02 57M27)},
	MRNUMBER = {3764859},
	MRREVIEWER = {Leandro\ Vendramin},
	DOI = {10.1080/00927872.2017.1316857},
	URL = {https://doi.org/10.1080/00927872.2017.1316857},
}

@book {dummit,
	AUTHOR = {Dummit, David S. and Foote, Richard M.},
	TITLE = {Abstract algebra},
	EDITION = {Third},
	PUBLISHER = {John Wiley \& Sons, Inc., Hoboken, NJ},
	YEAR = {2004},
	PAGES = {xii+932},
	ISBN = {0-471-43334-9},
	MRCLASS = {00-01 (16-01 20-01)},
	MRNUMBER = {2286236},
}

@article {dihedral,
	AUTHOR = {Elhamdadi, Mohamed and Macquarrie, Jennifer and Restrepo,
	Ricardo},
	TITLE = {Automorphism groups of quandles},
	JOURNAL = {J. Algebra Appl.},
	FJOURNAL = {Journal of Algebra and its Applications},
	VOLUME = {11},
	YEAR = {2012},
	NUMBER = {1},
	PAGES = {1250008, 9},
	ISSN = {0219-4988,1793-6829},
	MRCLASS = {20B25 (20N02)},
	MRNUMBER = {2900878},
	MRREVIEWER = {Alexander\ Zvonkin},
	DOI = {10.1142/S0219498812500089},
	URL = {https://doi.org/10.1142/S0219498812500089},
}

@article {centerless,
	AUTHOR = {Bardakov, Valeriy G. and Nasybullov, Timur and Singh, Mahender},
	TITLE = {Automorphism groups of quandles and related groups},
	JOURNAL = {Monatsh. Math.},
	FJOURNAL = {Monatshefte f\"ur Mathematik},
	VOLUME = {189},
	YEAR = {2019},
	NUMBER = {1},
	PAGES = {1--21},
	ISSN = {0026-9255,1436-5081},
	MRCLASS = {20N02 (20B25 57M27)},
	MRNUMBER = {3948284},
	MRREVIEWER = {David\ Stanovsk\'y},
	DOI = {10.1007/s00605-018-1202-y},
	URL = {https://doi.org/10.1007/s00605-018-1202-y},
}

@book {robinson,
	AUTHOR = {Robinson, Derek John Scott},
	TITLE = {A course in the theory of groups},
	SERIES = {Graduate Texts in Mathematics},
	VOLUME = {80},
	PUBLISHER = {Springer-Verlag, New York-Berlin},
	YEAR = {1982},
	PAGES = {xvii+481},
	ISBN = {0-387-90600-2},
	MRCLASS = {20-01},
	MRNUMBER = {648604},
}

@article {dual,
	AUTHOR = {Burrows, Wayne and Tuffley, Christopher},
	TITLE = {The rack congruence condition and half congruences in racks},
	JOURNAL = {J. Knot Theory Ramifications},
	FJOURNAL = {Journal of Knot Theory and its Ramifications},
	VOLUME = {35},
	YEAR = {2026},
	NUMBER = {4},
	PAGES = {Paper No. 2550086},
	ISSN = {0218-2165,1793-6527},
	MRCLASS = {20N02},
	MRNUMBER = {5036133},
	DOI = {10.1142/S0218216525500865},
}

@article {free-z,
	AUTHOR = {Farinati, Marco Andrés and Guccione, Jorge A. and Guccione, Juan J.},
	TITLE = {The homology of free racks and quandles},
	JOURNAL = {Comm. Algebra},
	FJOURNAL = {Communications in Algebra},
	VOLUME = {42},
	YEAR = {2014},
	NUMBER = {8},
	PAGES = {3593--3606},
	ISSN = {0092-7872,1532-4125},
	MRCLASS = {20N99 (55N35)},
	MRNUMBER = {3196064},
	MRREVIEWER = {Mahender\ Singh},
	DOI = {10.1080/00927872.2013.790392},
	URL = {https://doi.org/10.1080/00927872.2013.790392},
}

@misc{my-code,
	author   	= {L{\d\uhorn}c Ta},
	title    	= {\texttt{GL-Rack-Classification}},
	year     	= {2025},
	url      	= {https://github.com/luc-ta/GL-Rack-Classification},
	note        = {\url{https://github.com/luc-ta/GL-Rack-Classification}. Accessed: 2025-7-17.}
}

@article {takasaki-aut,
	AUTHOR = {Bardakov, Valeriy G. and Dey, Pinka and Singh, Mahender},
	TITLE = {Automorphism groups of quandles arising from groups},
	JOURNAL = {Monatsh. Math.},
	FJOURNAL = {Monatshefte f\"ur Mathematik},
	VOLUME = {184},
	YEAR = {2017},
	NUMBER = {4},
	PAGES = {519--530},
	ISSN = {0026-9255,1436-5081},
	MRCLASS = {22E40 (20B25 20N02 57M27)},
	MRNUMBER = {3718201},
	MRREVIEWER = {Osman\ Mucuk},
	DOI = {10.1007/s00605-016-0994-x},
	URL = {https://doi.org/10.1007/s00605-016-0994-x},
}

@article {new-atlas,
	AUTHOR = {Petkova, Ina and Schwartz, Noah},
	TITLE = {A {L}egendrian knot atlas for knots of arc index 10},
	JOURNAL = {Exp. Math.},
	FJOURNAL = {Experimental Mathematics},
	VOLUME = {35},
	YEAR = {2026},
	NUMBER = {1},
	PAGES = {111--259},
	ISSN = {1058-6458,1944-950X},
	MRCLASS = {57K33},
	MRNUMBER = {5046357},
	DOI = {10.1080/10586458.2024.2430715},
}

@article {virtual,
	AUTHOR = {Cattabriga, Alessia and Nasybullov, Timur},
	TITLE = {Virtual quandle for links in lens spaces},
	JOURNAL = {Rev. R. Acad. Cienc. Exactas F\'is. Nat. Ser. A Mat. RACSAM},
	FJOURNAL = {Revista de la Real Academia de Ciencias Exactas, F\'isicas y
	Naturales. Serie A. Matematicas. RACSAM},
	VOLUME = {112},
	YEAR = {2018},
	NUMBER = {3},
	PAGES = {657--669},
	ISSN = {1578-7303,1579-1505},
	MRCLASS = {57M27 (08A99 20N02)},
	MRNUMBER = {3819722},
	MRREVIEWER = {Jie\ Wu},
	DOI = {10.1007/s13398-017-0445-0},
	URL = {https://doi.org/10.1007/s13398-017-0445-0},
}

@article {virtual-leg,
	AUTHOR = {Cahn, Patricia and Levi, Asa},
	TITLE = {Vassiliev invariants of virtual {L}egendrian knots},
	JOURNAL = {Pacific J. Math.},
	FJOURNAL = {Pacific Journal of Mathematics},
	VOLUME = {273},
	YEAR = {2015},
	NUMBER = {1},
	PAGES = {21--46},
	ISSN = {0030-8730,1945-5844},
	MRCLASS = {57M27},
	MRNUMBER = {3290443},
	MRREVIEWER = {Micah\ Whitney\ Chrisman},
	DOI = {10.2140/pjm.2015.273.21},
	URL = {https://doi.org/10.2140/pjm.2015.273.21},
}

@article {magmoidal,
	AUTHOR = {Davydov, Alexei},
	TITLE = {Nuclei of categories with tensor products},
	JOURNAL = {Theory Appl. Categ.},
	FJOURNAL = {Theory and Applications of Categories},
	VOLUME = {18},
	YEAR = {2007},
	PAGES = {No. 16, 440--472},
	ISSN = {1201-561X},
	MRCLASS = {18D10},
	MRNUMBER = {2369108},
	MRREVIEWER = {Eric\ C.\ Rowell},
}

@article {matveev,
	AUTHOR = {Matveev, S. Vladimir},
	TITLE = {Distributive groupoids in knot theory},
	JOURNAL = {Mat. Sb. (N.S.)},
	FJOURNAL = {Matematicheski\u i\ Sbornik. Novaya Seriya},
	VOLUME = {119(161)},
	YEAR = {1982},
	NUMBER = {1},
	PAGES = {78--88, 160},
	ISSN = {0368-8666},
	MRCLASS = {57M25 (20L15)},
	MRNUMBER = {672410},
	MRREVIEWER = {Jonathan\ A.\ Hillman},
}

@misc{ta,
	title={Good involutions of conjugation subquandles}, 
	author={Ta, L{\d\uhorn}c},
	year={2025},
	eprint={2505.08090},
	archivePrefix={arXiv},
	primaryClass={math.GT},
	url={https://arxiv.org/abs/2505.08090}, 
	NOTE = {Preprint, arXiv:2505.08090 [math.GT].},
}

@article {tensors,
	AUTHOR = {Kamada, Seiichi},
	TITLE = {Tensor products of quandles and 1-handles attached to
	surface-links},
	JOURNAL = {Topology Appl.},
	FJOURNAL = {Topology and its Applications},
	VOLUME = {301},
	YEAR = {2021},
	PAGES = {Paper No. 107520, 18},
	ISSN = {0166-8641,1879-3207},
	MRCLASS = {57K12},
	MRNUMBER = {4312970},
	MRREVIEWER = {Indu\ Rasika\ Churchill},
	DOI = {10.1016/j.topol.2020.107520},
	URL = {https://doi.org/10.1016/j.topol.2020.107520},
}

@misc{wood,
	title={Distinguishing power of 4-{L}egendrian permutation racks}, 
	author={Ta, L{\d\uhorn}c and Wood, Peyton},
	year={2026},
	eprint={2510.26619},
	archivePrefix={arXiv},
	primaryClass={math.GT},
	url={https://arxiv.org/abs/2510.26619}, 
	NOTE = {Preprint, arXiv:2510.26619 [math.GT].},
}

@article {cheng,
	AUTHOR = {Cheng, Zhiyun and He, Zhiyi},
	TITLE = {Fundamental generalized {L}egendrian rack and classical
	invariants},
	JOURNAL = {Internat. J. Math.},
	FJOURNAL = {International Journal of Mathematics},
	VOLUME = {37},
	YEAR = {2026},
	NUMBER = {4},
	PAGES = {Paper No. 2650027, 14},
	ISSN = {0129-167X,1793-6519},
	MRCLASS = {57K12},
	MRNUMBER = {5051216},
	DOI = {10.1142/S0129167X26500278},
}

@misc{cheng2,
	title={{GL}-racks and coloring invariants of {L}egendrian knots}, 
	author={Zhiyun Cheng and Zhiyi He},
	year={2026},
	eprint={2605.16968},
	archivePrefix={arXiv},
	primaryClass={math.GT},
	url={https://arxiv.org/abs/2605.16968}, 
	NOTE = {Preprint, arXiv:2605.16968 [math.GT].},
}
		
	\end{document}